\documentclass{amsart}
\usepackage{amssymb}
\usepackage{amsfonts}

\theoremstyle{plain}

\newtheorem{corollary}{Corollary}

\newtheorem{definition}{Definition}

\newtheorem{lemma}{Lemma}

\newtheorem{proposition}{Proposition}
\newtheorem{remark}{Remark}

\numberwithin{equation}{section}
\newcommand {\fp}{\mathfrak{p}}
\newcommand {\fg}{\mathfrak{g}}
\newcommand {\fk}{\mathfrak{k}}
\newcommand {\ep}{\epsilon}
\newcommand {\IR}{\mathbb{R}}
\newcommand {\IH}{\mathbb{H}}
\newcommand{\IS}{\mathbb{S}}
\begin{document}
\title[ Integrable systems]{Integrable Hamiltonian systems on symmetric
spaces:\\
Jacobi, Kepler and Moser}
\author{ Velimir Jurdjevic}
\address{ \ Department of Mathematics, University of Toronto}
\email{jurdj@math.toronto.edu}
\curraddr{Department of Mathematics, University of Toronto, 40 St. George
st.. Toronto}
\urladdr{http://www.jurdj.math.toronto.edu}
\date{October 15, 2009}
\dedicatory{Dedicated to the memory of J. Moser.}
\thanks{This is a working copy of the paper that will be submitted for a
publication upon its completion. }

\begin{abstract}
This paper defines a class of left invariant variational problems on a Lie
group$\ G\ $whose Lie algebra $\mathfrak{g}$\ admits Cartan decomposition $%
\mathfrak{g}=\mathfrak{p}+\mathfrak{k}$\ with the usual Lie algebraic
conditions 
\begin{equation*}
\lbrack \mathfrak{p},\mathfrak{p]\subseteq k\ },\ \mathfrak{[p},\mathfrak{%
k]\subseteq p},\mathfrak{\ [k},\mathfrak{k]\subseteq k.}
\end{equation*}

The Maximum Principle of optimal control leads to the Hamiltonians \ $H$ on $%
\mathfrak{g\ }$that admit spectral parameter representations with \
important contributions to the theory of integrable Hamiltonian systems. \
Particular cases provide natural explanations
for the classical results of Fock and Moser linking Kepler's problem to the
geodesics on spaces of constant curvature and \ J.Moser's work on integrability based on isospectral
methods in which C. Newmann's  mechanical problem on the sphere and C. L. Jacobi's geodesic problem on an ellipsoid  play the central role. \ The paper also shows the relevance of this class of Hamiltonians
to the elastic curves on \ spaces of constant curvature.
\end{abstract}

\maketitle

\section{\protect\bigskip Introduction}

$\ $ \ A Lie group $G\ $with an involutive automorphism$\ \sigma \ $admits\
several$\ $natural\ variational\ problems\ whose\ solutions\ provide\ new\
insights\ into the theory of integrable Hamiltonian systems and to the
geometry of the associated homogeneous spaces.\ An \ involutive automorphism 
$\sigma $\ on a Lie group $G\ $induces a splitting $\mathfrak{g=p\oplus k\ }$%
of the Lie algebra $\mathfrak{g\ }$of$\ G\mathfrak{\ }$with $\mathfrak{k\ }$%
equal to the Lie algebra of the group $K\ $of fixed points under $\sigma .\ $%
When $G\ $\ is semisimple then $\mathfrak{p\ }$is the orthogonal complement
to $\mathfrak{k\ }$ \ relative to the Killing form  $\mathfrak{\ }$and \ $%
\mathfrak{p}$\ and $\mathfrak{k\ }$\ satisfy the following Lie algebraic
relations:%
\begin{equation}
\lbrack \mathfrak{p},\mathfrak{p]=k\ },\ \mathfrak{[p},\mathfrak{k]=p.}
\end{equation}%
The \ first relation implies that any two points of $G$ can be connected by a curve whose tangent  takes values  in the left invariant
distribution $\mathcal{D(}g)=\{gU:U\in \mathfrak{p\}.\ }$  In the case that $%
(G,K)\ $is a Riemannian symmetric pair there is an $Ad_{K}\ $invariant,
positive definite quadratic $\ $form $\langle \ ,\ \rangle \ $on$\ \mathfrak{p\ }$that \
induces a natural optimal control problem on $G:\ $ minimize the integral$\
\ \frac{1}{2}\int_{0}^{T}\langle U(t)\ ,U(t)\ \rangle dt\ \ $among all curves $g(t)\in G\ $
that are the solutions of%
\begin{equation}
\frac{dg}{dt}=g(t)U(t),\ U(t)\in \mathfrak{p},\ t\in \lbrack 0,T]
\label{subr}
\end{equation}%
with fixed \ boundary conditions $g(0)=g_{0}\ ,g(T)=g_{1}.\ $Here$\ 
\mathfrak{\ }g_{0}\ $and$\ \ g_{1}\ $are arbitrary but fixed points\ in $G\
\ $and the terminal time $T>0\ $\ is also fixed. $\ $\ This problem, \
  called \ the\ \emph{\ }canonical sub-Riemannian problem \emph{\ }on$\emph{%
\ }G\emph{,}$\ is well defined in the sense that for any pair of boundary points in $G$ there is an optimal solution. It\ is well\ known\ that\  optimal  solutions are of the form\ 
\begin{equation}
g(t)=g_{0}e^{t(P+Q)}e^{-tQ}
\end{equation}%
for some elements $P\in \mathfrak{p}$,\ and\ $Q\in \mathfrak{k\ }$(\cite{Mem}%
\ ).  The above implies that any element $g\ $in $G\ $can be
represented as \ $g=e^{(P+Q)}e^{-Q}\ $for some elements $P\in \mathfrak{p}$%
,\ and\ $Q\in \mathfrak{k}.\ $This sub-Riemannian problem \ is naturally
related to the canonical Riemannian problem \ on the quotient space $M=G/K$\
in the sense that the Riemannian geodesics are the projections of the above
curves with $Q=0.$

In this paper we will be interested in another\ optimal control problem
defined \ by an affine distribution $\mathcal{D}(g)=\{g(A+U):U\in \mathfrak{%
k\}\ }$with $A\ $ a\ regular element $\ $in\ $\mathfrak{p\ }$under the
assumption that the Killing form is  definite on $\mathfrak{k\ }(\ $%
which is true when $K$\ is a compact subgroup\ of $G\ ).$ This\ optimal
control problem consists of finding the minimum of \ $\frac{1}{2}%
\int_{0}^{T}\langle U(t),U(t)\rangle dt\ \ $among all \ solution curves $g(t)\in G\ $ of\
the affine control problem

$\ \ $%
\begin{equation}
\frac{dg}{dt}(t)=g(t)(A+U(t)),U(t)\in \mathfrak{k},\ t\in \lbrack 0,T]~
\label{Affine}
\end{equation}

subject to the given boundary conditions $g(0)=g_{0},\ g(T)=g_{1}\ $where
the quadratic form $\langle\,,\,\rangle$ denotes a  scalar multiple of the Killing form that is positive definite on $\fk$.

This problem might be regarded as the\ canonical$\ $affine\ problem on
symmetric pairs $(G,K)$ for the following reasons.  The reciprocal affine system with $A\in 
\mathfrak{k}$ and\ $U\in \mathfrak{p}$\  is isomorphic  to (\ref{subr}%
) and bears little resemblance to the solutions of (\ref{Affine}).  Moreover, \ two affine systems\ \ defined by regular elements\ $A_{1}$\ and $%
A_{2}\ $are conjugate.\ 

The\ affine problem (\ref{Affine}) will be referred to
as the Affine-Killing problem or \ (\textbf{Aff}) for brevity. 
It is first shown that\ (\textbf{Aff}) is well defined in the sense that for
any pair of boundary conditions$\ (g_{0},g_{1})\ $there\ exist\ $T>0\ $and a
solution $g(t)~$of\ (\ref{Affine}) that satisfies$\ g(0)=g_{0},\ g(T)=g_{1}\ 
$such\ that\ the the control $U(t)\ $ that generates $g(t)$ minimizes the integral \ $\frac{1}{2}%
\int_{0}^{T}\langle U(t),U(t)\rangle dt\ $among all other controls whose solution curves  satisfy the
same boundary data. Then it \ is shown that optimal solutions are the
projections of the integral curves of a certain Hamiltonian system on the
cotangent bundle\ \ $T^{\ast }G$\ of $G$\ obtained through the use of the
Maximum Principle of Optimal Control. 

To preserve the left invariant
symmetries, $T^{\ast }G$\ \ is realized as the product $%
G\times \mathfrak{g}^{\ast }$\ with \ $\mathfrak{g}^{\ast }\ $equal to the
dual of the Lie algebra \ $\mathfrak{g\ }$of $G$ and then $\mathfrak{g}%
^{\ast }$ is identified with \ $\mathfrak{g\ }$via the\ Killing\ form so
that ultimately\ $T^{\ast }G\ $is realized as $G\times \mathfrak{g}$. In
this setting the Hamiltonian associated with \ (\textbf{Aff) }is of the
form\ \ 

\begin{equation*}
H=\frac{1}{2}\langle L_{\mathfrak{k}},L_\fk\rangle+\langle A,L_{\mathfrak{p}}\rangle
\end{equation*}%
\ \ \ 

where $L_{\mathfrak{p}}$ and \ $L_{\mathfrak{k}}\ $\ denote the projections
of $L\in \mathfrak{g\ }$onto the factors $\mathfrak{p}$\ and $\mathfrak{k}$%
.\  

\ \ \ Because $K\ $acts on $\mathfrak{p\ }$via the adjoint action,\ \ $%
\mathfrak{g}$ \ as a vector space carries two Lie algebras: a Lie algebra of 
$G$\ and the Lie algebra $\mathfrak{g}_{s}\ $of the semidirect product $%
G_{s}=K\rtimes \mathfrak{p.\ }$\  \ The \ \ affine problem\emph{\ }%
 then admits an analogous formulation \ on $G_{s}$  and  the
Maximum Principle leads to the Hamiltonian $H$\ that formally looks the same as
the one obtained in the semisimple case. \ We refer to the semidirect version of\ (%
\textbf{Aff)} \ as the\emph{\ semidirect shadow problem}.\ 

The essence of the paper lies in the integrability properties of the
associated Hamiltonian flows\ $\vec{H}$\ and$\ \vec{H}_{s}.\ $ It is shown
that \ each flow admits a spectral representation 
\begin{equation}
\frac{dL_{\lambda }}{dt}=[M_{\lambda },L_{\lambda }]
\end{equation}

where \ $M_{\lambda }=\frac{1}{\lambda }(L_{\mathfrak{p}}-\epsilon A)\ $and\ 
$L_{\lambda }=L_{\mathfrak{p}}-\lambda L_{\mathfrak{k}}+(\lambda
^{2}-s )A,$\ with $s =1\ $in the semisimple case and $s
=0\ $in the semidirect case. Hence, the spectral invariants are 
constants of motion for the associated Hamiltonian flows.

On spaces of constant curvature these results recover the integrability
results \ associated with elastic curves and their \ mechanical counterparts
\ ( \cite{Mem}).\ 
Remarkably,\ the spectral invariants above also recover the classical\
integrability \ results C. Newmann \ for mechanical systems \ with quadratic
potential (\cite{Nmn})\ and the related results of C.G.J. Jacobi \
concerning the geodesics on an ellipsoid. The present formalism also
clarifies the contributions of J. Moser ( \cite{MosCh})\ \ on integrability
of Hamiltonian systems based on isospectral methods. More significantly,
this study reveals a large class of integrable Hamiltonian systems in which
these classical examples appear only as very particular cases. 

\section{Notations and the Background material}

The basic setting is most naturally defined through the language of
symmetric spaces.  The  essential ingredients are assembled below.

\ An
involutive automorphism $\sigma \ $on $G\ $is an analytic mapping $%
G\rightarrow G,\ \sigma \neq I$\ that satisfies%
\begin{equation}
\ \sigma (g_{2}g_{1})=\sigma (g_{2})\sigma (g_{1}),\ \ for\ all\
g_{1},g_{2}\ in\ G.
\end{equation}

Then the tangent map $\mathfrak{\sigma }_{\ast }\ $of $\sigma $ induces a
splitting $\mathfrak{g=p\oplus k\ }$of the Lie algebra\ $\mathfrak{g}$ of $G$
\ with \ 
\begin{equation}
\mathfrak{p=}\{A\in \mathfrak{g:\sigma }_{\ast }(A)=-A\}\ and\ k\mathfrak{=\{%
}A\in \mathfrak{g:\sigma }_{\ast }(A)=A\}
\end{equation}

\ The fact that $\mathfrak{\sigma }_{\ast }$ is\ a Lie algebra automorphism
\ easily implies the following Lie algebraic relations%
\begin{equation}
\lbrack \mathfrak{p},\mathfrak{p]\subseteq k\ },\ \mathfrak{[p},\mathfrak{%
k]\subseteq p},\mathfrak{\ [k},\mathfrak{k]\subseteq k}  \label{cartan}
\end{equation}

It follows that $\mathfrak{k\ }$is a Lie subalgebra of $\mathfrak{g,\ }$%
equal to the Lie algebra of the group $K=\{g\in G:\sigma (g)=g\}\ $and that $%
\mathfrak{p}$\ is an\ $Ad_{K}$ \ invariant\ vector subspace of $\mathfrak{g}$%
 \ in the sense that $Ad_{h}(\mathfrak{p})\subseteq \ \mathfrak{p}$\ for
any $h\in K.\ $An$\ Ad_{K}$ \ invariant non-degenerate $\ $quadratic form on 
$\mathfrak{p\ }$will be called \emph{pseudo Riemannian.}\ It is easy to
show\ by differentiating that an $Ad_{K}$ invariant \ quadratic form $<\ ,\
>\ $on\ $\ \mathfrak{p\ }$is \emph{invariant}, in the sense that%
\begin{equation}
<A,[B,C]>=<[A,B],C>
\end{equation}%
for any $A,C\ $in$\ \mathfrak{p}$\ and any $B$\ in$\ \mathfrak{k}.$

A pseudo Riemannian form that is positive definite will be called \emph{%
Riemannian.\ }\ In the literature of symmetric spaces (\cite{Hel}) the
pair $(G,K)$, \ with $K$ a closed subgroup of $G\ $\ obtained by an
involutive automorphism on $G\ \ $described above, is called a \emph{%
symmetric pair.}\ If in addition \ this pair admits an $Ad_{K}$ invariant \
positive definite quadratic form $\langle \ ,\ \rangle \ $on\ $\mathfrak{p\ \ }\ $then it
is called\ a \emph{Riemannian symmetric pair.}\ \ \ Riemannian symmetric
pairs can be characterized as follows.
\begin{proposition} Let\ $Ad_{h,\mathfrak{p}}$ denote the
restriction of $Ad_{h}\ $to $\mathfrak{p.}$\ Then \ a symmetric pair $(G,K)$%
\ admits a Riemannian \ quadratic form $\langle \ ,\ \rangle\ $on\ $\ \mathfrak{p\ }$\ if
and only if \\$\{Ad_{h,\mathfrak{p}}:\ h\in K\}$ is a
\ compact subgroup of $Gl(\mathfrak{p}).$\ 
\end{proposition}

In the text below we will make use of the Killing form$\
\langle A,B\rangle_{k}=Tr(adA\circ adB)$\ for $A\ $and\ $B\ $in\ $\mathfrak{g}$,\
where $Tr(X)$ denotes the trace of a linear endomorphism $X.$
The Killing form is \ invariant under any automorphism $\phi $\ on\ $%
\mathfrak{g}\ $\ and in particular it is\ $Ad_{K}$ and\ $Ad_{G}~$invariant\ (%
\cite{Hel}).\ The invariance relative to $Ad_{G}$ implies
 \begin{equation}%
\langle A,[B,C]\rangle _{k}=\langle [A,B],C]\rangle _{k}\end{equation}\ for any matrices $A,B,C\ $in\ $\mathfrak{g}.\ 
$\ 
Spaces$\ \mathfrak{p}$\ and\ $\mathfrak{k\ }$are
orthogonal relative to $\langle \ ,\ \rangle _{k}$ because$\ \langle A,B\rangle _{k}=\langle \sigma _{\ast
}(A),\sigma _{\ast }(B)\rangle _{k}=\langle -A,B\rangle_{k}=-$\ $\langle A,B\rangle _{k}\ $for any $A\ $in$\ 
\mathfrak{p}$\ and any $B$\ in$\ \mathfrak{k}.\ $

In this paper $(G,K)$\ will be assumed  a symmetric Riemannian pair
with $G$\ semisimple and connected$\ $\ and \ $K\ $compact. \
Semisimplicity implies that the Killing form is non-degenerate, \ which then
implies that its restriction to $\mathfrak{p\ }\ $is pseudo Riemannian. \
Semisimplicity also implies that the Cartan relations (\ref{cartan}) take on a
stronger form 
\begin{equation}
\lbrack \mathfrak{p},\mathfrak{p]=k\ },\ \mathfrak{[p},\mathfrak{k]=p},%
\mathfrak{\ [k},\mathfrak{k]\subseteq k.}  \label{cartan1}
\end{equation}

The fact that $K\ $is a compact subgroup of $G\ $implies that the Killing
form is negative definite on $\mathfrak{k}$ (\cite{Eb}, p. 56). In the sequel $%
\langle\ ,\ \rangle\ $will\ denote$\ $any scalar multiple of the Killing form which is
positive definite on $\mathfrak{k.}$ Under these conditions then $||U||$ will denote the induced norm $||U||=\sqrt{\langle U,U\rangle}$.

An element $A$ in $\mathfrak{p\ }$is
said to be\emph{\ regular} if $\{B\in \mathfrak{p}:[A,B]=0\}$\ is an abelian
algebra. It follows that $A$\ is regular if and only if the algebra $\mathbb{A}$ spanned by $\{B\in 
\mathfrak{p}:[A,B]=0\}\ $is a maximal abelian algebra in $\mathfrak{p}$ that
contains $A$\ (\cite{Eb}\ ).

With these notions at our disposal we $\ $return \ now to the\ affine problem%
\emph{\ }defined above.$\ $\ It will be convenient to adopt the language of
control theory and \ regard \ (\ref{Affine})\ as a control system with $%
U(t)\ $\ playing the role of \emph{control}. \ In order to meet the
conditions of the Maximum Principle control functions are assumed bounded and measurable on compact intervals $[0,T]$. Solutions of (\ \ref{Affine} ) are called \emph{%
trajectories. }\ A control $U(t)\ $\ is said to \emph{steer\ }$g_{0}\ $to $%
g_{1}\ $in $T\ $units of time if the corresponding trajectory \ $g(t),t\in
\lbrack 0,T]$\ \ satisfies $g(0)=g_{0},\ g(T)=g_{1}.$ A trajectory\ $g(t)\ $ generated by a control $U(t)$ 
 on an interval $[0,T]\ $ is \emph{optimal\ relative to the boundary
conditions (g}$_{0},g_{1})\ $\emph{\ }if the  integral $\frac{1}{2}$ $\int_{0}^{T}||U(t)||^2\,dt\ $is minimal
among all other controls that steer\emph{\ }$g_{0}\ $to $g_{1}\ $in $T\ $%
units of time. \ Controls that result in optimal trajectories \ are called 
\emph{optimal.\ }Thus every optimal control $U(t)\ $gives rise to a unique
optimal trajectory because  the initial point $g_0$ is fixed.
\section{The existence of optimal solutions}

\begin{proposition}
\label{Contr} If $\ A$\ is regular\ then (\textbf{Aff)} problem is well
posed $\ $in the sense that\ for any pair of boundary conditions \ $%
(g_{0},g_{1})\ $there exist a time $T>0$\ and a solution\ $g(t)\ $on the
interval $[0,T]\ $that is optimal relative to $g_{0}\ $and\ $g_{1}.\ $
\end{proposition}

The proof of this proposition \ requires   several auxiliary facts  from the optimal control theory and from the theory of symmetric spaces. We begin  first with the facts from the theory of symmetric spaces (\cite{Eb},\ \cite{Hel}). 

A  Lie algebra $\mathfrak{g}\ $ is said to be simple if it contains no ideals
other that $\{0\}\ $and\ $\mathfrak{g}$.\ A Lie group $G$\ is said to be
simple if its Lie algebra is simple.\ $\ $The first fact is given by

\begin{lemma}
\label{irred} If $(G,K)$ is a symmetric pair with $G$ simple  then $Ad_{K}$\ acts irreducibly on $\mathfrak{p}\ $.
\end{lemma}

\begin{proof}
Let $V$\ denote an $Ad_{K}\ $invariant vector subspace of $\mathfrak{p}.\ $%
Denote by $V^{\perp }$ the orthogonal complement of $V$\ in $\mathfrak{p\ }$%
relative to the Killing form\ $\langle \ ,\ \rangle_{k}.$\ Since

$\langle V,Ad_{h}(V^{\perp })\rangle _{k}=$ $\langle Ad_{h^{-1}}(V),V^{\perp }\rangle _{k},$\ it follows
that $V^{\perp }\ $is also $Ad_{K}\ $invariant. Therefore, $[\mathfrak{k}%
,V]\subseteq V$\ and $[\mathfrak{k},V^{\perp }]\subseteq V^{\perp },\ $which
in turn implies that $\langle\mathfrak{k},[V,V^{\perp }]\rangle _{k}=0.$ It follows that $%
[V,V^{\perp }]=0\ $by semisimplicity of $\mathfrak{g}.$

The above implies that $V+[V,V]$\ is an ideal in $\mathfrak{g\ }$that is
orthogonal to $V^{\perp }.\ $Since $\mathfrak{g\ }$is simple\ $\ V+[V,V]=%
\mathfrak{g}$\ and$\ $\ therefore $V^{\perp }=0.\ $But then $V=\mathfrak{p}$.
\end{proof}

The other facts which are needed for the proof are assembled below in the
forms of propositions.

\begin{proposition}
\label{decomp1} Suppose that $(G,K)$\ is\ a\ symmetric Riemannain pair.
There exist linear subspaces $\ \mathfrak{p}_{1},$ $\mathfrak{p}_{2},\ldots
, $ $\mathfrak{p}_{m}\ $\ of $\fp$ such that

(1) $\mathfrak{p=p}_{1}\oplus $ $\mathfrak{p}_{2}\oplus \cdots \cdots \oplus 
$ $\mathfrak{p}_{m}$ and

(2) $\mathfrak{p}_{1},$ $\mathfrak{p}_{2},\ldots ,$ $\mathfrak{p}_{m}\ $are
pairwise orthogonal relative to the Killing form$.$

(3) Each $\mathfrak{p}_{i}$\ is $ad(\mathfrak{k)}\ $invariant and contains
no proper $ad(\mathfrak{k)}$\ invariant \ linear subspace.
\end{proposition}

\begin{proposition}
\label{decomp2} Let $\mathfrak{g}_{i}=%
\mathfrak{p}_{i}$\ \ $+[$ $\mathfrak{p}_{i}$\ $,$ $\mathfrak{p}_{i}$\ $]\ $%
 $i=1,\ldots ,m \ $ in a semisimple Lie algebra $\fg$. Then,
(1) Each $\mathfrak{g}_{i}\ $is an ideal of $\mathfrak{g}\ $\  and also a
simple Lie\ algebra.\ 

Moreover,

(2) $[$ $\mathfrak{g}_{i},$ $\mathfrak{g}_{j}]=0\ \ $and \ $\langle$ $\mathfrak{g}%
_{i},$ $\mathfrak{g}_{j}\rangle _{k}=0,i\neq j$, and 
$\{X\in \mathfrak{g}:[X,Y]=0 \,,Y\in \mathfrak{g\}=0.}$
\end{proposition}

\begin{corollary}
Let $A_{i}$\ denote the projection of
a regular element $A\ $in$\ \mathfrak{p\ }$on $\fp_i$. Then each $A_{i}\neq 0, i=1,\dots,m.$
\end{corollary}

\begin{proof}
If $A_{i}$ were equal to $0\ $then $[A,\mathfrak{p}_{i}]=0\ $by\ (2) in
Proposition \ref{decomp1}. \ This would imply that $\mathfrak{p}_{i}\ $is
abelian by regularity of $A$, which in turn would imply that $\mathfrak{g}%
_{i}=$ $\mathfrak{p}_{i}.$\ But that would contradict (3) in Proposition \ref%
{decomp1}.
\end{proof}

We now turn attention to the pertinent ingredients from the
accessibility \ theory of control systems. \ The Lie Saturate of a left invariant family of
vector fields $\mathcal{F\ }$ is the largest family of left invariant vector fields (in the sense of set inclusion ) that leaves the closure of the reachable sets of $\mathcal{F}$ invariant (\cite{jurd}). It is
denoted by $LS(\mathcal{F}).$ 

Since  left invariant vector fields are defined by their values at the identity,  the Lie saturate  admits a paraphrase in terms of  the defining  set $\Gamma$ in $\fg$.  For  the affine system (\ref{Affine}\ ),
 $\Gamma =\{A+B:B\in \mathfrak{k}\}.$
\begin{definition} Let $\Gamma\subseteq\fg$ . The reachable set of $\Gamma$ denoted by $A(\Gamma)$  is the set of terminal points $ g(T)\in G$  corresponding to
the absolutely continuous curves $g(t)$ on  intervals $[0,T]$ such that $ g(0)=I$
and $\frac{dg}{dt}(t)g^{-1}(t)\in \Gamma$  for almost all $t\in [ 0,T].$
\end{definition}

Then \ $LS(\Gamma ),$\
the Lie Saturate of $\ \Gamma \ $  can be described as the largest family in $\mathfrak{g}$\
such that \ 
\begin{equation}
cl(\mathcal{A}(LS(\Gamma ))=cl(\mathcal{A}(\Gamma ),
\end{equation}
where $cl(X)$ denotes the topological closure of a set $X$.

The following lemma is well known in control theory (\cite{jurd}).

\begin{lemma}
\label{sat}\ a.\ $LS(\Gamma )=\mathfrak{g}$\ is a necessary and sufficient
condition that $\mathcal{A}(\Gamma )=G.$\ 

If $C$ denotes the convex cone spanned by $%
\sum_{i=1}^m \alpha _{i}Ad_{h_{i}}(A),\ h_{i}\in \mathfrak{k},\ \alpha _{i}\geqq
0,i=1,2,\ldots ,m, m\in Z^{+} $ then 
\begin{equation}
C\cup \mathfrak{k}\subseteq LS(\Gamma )
 \end{equation}where 
  $\Gamma=\{A+U:U\in \fk\}.$

\end{lemma}

With these results at our disposal let us turn to the proof of the
proposition.

\begin{proof}
Let $Traj(g_{0},g_{1})$ denote the set of solutions $g(t)$ of \ (\ref{Affine}%
) that satisfy$\ g(0)=g_{0},\ g(T)=g_{1}$\ for some $T>0.$ If \ $%
Traj(g_{0},g_{1})$ is not empty for any $g_{0}\ $and$\ g_{1}$ in $G\ $then (%
\ref{Affine}) is said to be \emph{controllable}.  An argument based on  weak compactness of  closed balls in Hilbert spaces 
shows that there is an optimal \
trajectory $\hat{g}(t)\ $in\ $Traj(g_{0},g_{1})\ $generated by a control
 $\hat{U}(t)\ $in $L^{2}([0,T])$\ whenever\ 
$Traj(g_{0},g_{1})$\ is not empty (Theorem 1 in (\cite{Mem}). But then it can be shown that an optimal control in $%
L^{2}([0,T])$ is absolutely continuous and hence belongs to  $%
L^{\infty }([0,T]).\ $The above argument shows that  controllability implies the existence of optimal trajectories.

 Now address the question of controllability. According to Lemma 2 it would suffice to show that  the convex cone $C$ spanned by $\{\sum \alpha _{i}Ad_{h_{i}}(A),\ h_{i}\in \mathfrak{k},\ \alpha
_{i}\geqq 0\ \}\ $is equal to $\mathfrak{p }$.

Let $V$ denote the vector space spanned by $\{Ad_{h}(A):h\in \mathfrak{k\ }%
\}\ $and let $V_{i}\ $denote the vector space spanned by $%
\{Ad_{h}(A_{i}):h\in \mathfrak{k\ }\}\ $where $A_{i}$\ is the projection of $%
A$ on $\mathfrak{p}_{i}$, as in the Corollary above.\ Each $V_{i}$ is a
non-zero $Ad_{K}$\ \ invariant vector subspace of a simple Lie algebra $%
\mathfrak{g}_{i}.\ $According to Lemma \ref{irred}, $V_{i}=$ $\mathfrak{p}%
_{i}, i=1,\dots,m $ and hence, $V=\mathfrak{p.}$ It follows that \ \ $C=\{\sum \alpha
_{i}Ad_{h_{i}}(A),\ h_{i}\in \mathfrak{k},\ \alpha _{i}\geqq 0\ \}\ $is an$\
Ad_{K}$ \ invariant convex cone with a non empty interior in $\mathfrak{p.}$

\ Then $C=\mathfrak{p}$\ if and only if the origin in $\mathfrak{p\ }$were
contained in$\ $the interior of $C.$\ Let $S^{n}=\{X\in \mathfrak{p%
}:\Vert X\Vert =\Vert A\Vert \}.\ $If $0$\ were not in the the interior of $%
C,\ $then $C\cap S^{n}$ would be a convex cone in the sense of Eberlein (%
\cite{Eb}, 1.15)\ that is invariant under $Ad_{K}.\ $But then the sole of
this convex set would be a fixed point of $Ad_{K}\ $which is not possible
since $Ad_{K}$\ acts irreducibly on each $\mathfrak{p}_{i}.$
\end{proof}

\subsection{Semidirect products and the shadow\ problem}

\ Recall that if $K_{0}\ $is a Lie group which acts linearly on \ a finite
dimensional vector space $V\ $then \ the semidirect product $G_{s}=V\rtimes
K_{0}$\ \ consists of points ($v,h)\ $in \ $V\times K_{0}\ $with the group
operation \ $(v_{1},h_{1})(v_{2},h_{2})=(v_{1}+h_{1}(v_{2}),h_{1}h_{2}).$ The
Lie algebra $\mathfrak{g}_{s}\ $of $G_{s}\ $consists of pairs $(A,B)\ $in\ $%
V\times \mathfrak{k}_{0}\mathfrak{\ }$where $\mathfrak{k}_{0}\mathfrak{\ }$%
denotes the Lie algebra of $K_{0},\ $with the Lie bracket $%
[(A_{1},B_{1}),(A_{2},B_{2}]_{s}=(B_{1}(A_{2})-B_{2}(A_{1}),[B_{1},B_{2}]).\ 
$

Every$\ $ semidirect product $V\ltimes K_{0}$\ admits\ an\ involutive\
automorphism$\ \sigma (x,h)=(-x,h)\ $for every $(x,h)\in $ $V\ltimes K_{0}.$%
\ It follows that $K=\{0\}\times K_{0}\ $is the group of fixed points of $%
\sigma $ and that 
\begin{equation}
\mathfrak{p}=V\times \{0\}\ and\ \mathfrak{k}=\{0\}\times \mathfrak{{k}}%
_{0}.
\end{equation}
It is easy to check that $Ad_{h}(x,0)=(h(x),0)$\ for every $h\in K$\ and
every $x\in V.\ $Therefore,\ $\ (G,K)\ $is a symmetric Riemannian pair if
and only if $\ K_{0}\ $is a compact subgroup of $Gl(V).\ $

Every Lie group $G\ \ $that admits an~involutive automorphism~carries the
semidirect product $G_{s}=\mathfrak{p}\rtimes K\ $\ because \ $K\ $acts
linearly on the Cartan space $\mathfrak{p\ \ }$ via the transformation \ 
\begin{equation}
h(A)=Ad_{h}(A),A\in \mathfrak{p}\mbox{, for each }h\in K.
\end{equation} 
Therefore,  
the Lie bracket on $\mathfrak{g%
}_{s}$ is given by 
\begin{equation}
\lbrack
(A_{1},B_{1}),(A_{2},B_{2})]_s=(adB_{1}(A_{2})-adB_{2}(A_{1}),[B_{1},B_{2}]).
\end{equation}%
\ 
If  $(A,B)\ $in$\ \mathfrak{p}\rtimes \mathfrak{k\ }
$ is identified with $A+B\ $\ in \ $\mathfrak{p+k\ }$ then the semidirect Lie bracket $[\
,\ ]_{s}$\ can be \ redefined as 
\begin{equation}
\lbrack
(A_{1}+B_{1}),(A_{2}+B_{2})]_s=[B_{1},A_{2}]-[B_{2},A_{1}]+[B_{1},B_{2}],
\end{equation}

from which it follows that%
\begin{equation}
\lbrack \mathfrak{p},\ \mathfrak{p}]_{s}=0,\ [\mathfrak{p},\ \mathfrak{k}%
]_{s}=[\mathfrak{p},\ \mathfrak{k}],\ [\mathfrak{k},\ \mathfrak{k}]_{s}=[%
\mathfrak{k},\ \mathfrak{k}].
\end{equation}

Thus $\mathfrak{g}$\ is the underlying
 vector space for both Lie algebras\ $\mathfrak{g}$\ and\ $\mathfrak{g}_{s}$, a fact
which is important for the subsequent development.\ \ \ The passage from \ $%
\mathfrak{g_s}$\ \ to \ $\mathfrak{g}$\ \ can be described by a continuous
parameter $s $ by deforming the Lie algebra $\mathfrak{g_s\ }$to\ $%
\mathfrak{g}$\ via the Lie bracket\ $[\ ,\ ]_{s }\ :$%
\begin{equation}
\lbrack (A_{1}+B_{1}),(A_{2}+B_{2})]_{s
}=[B_{1},A_{2}]-[B_{2},A_{1}]+[B_{1},B_{2}]+ s\lbrack A_{1},A_{2}]
\label{defbr}
\end{equation}

We now return briefly to the affine problem\ {\textbf{Aff}} to note  that the data which is required for its formulation on a semisimple Lie group $G$ also permits  a formulation  on the semidirect product $G_{s}$. The semidirect version consists of 
minimizing the integral $\frac{1}{2}$ $\int_{0}^{T}\langle U(t),U(t)\rangle\,dt\ $ over all
solutions\ $g(t)\ $\ in $G_{s}\ $of

\begin{equation*}
\frac{dg}{dt}(t)=g(t)(A+U(t)),U(t)\in \mathfrak{k},\ t\in \lbrack 0,T]
\end{equation*}

that meet the boundary conditions $g(0)=g_{0},g(T)=g_{1}.$ This "shadow"\
 problem will be referred to  as
(\textbf{Aff}$_{s}$).$\ $The same arguments used in the semisimple case
show that (\textbf{Aff}$_{s})\ $is also well defined\ in the sense of
Proposition\ \ref{Contr}.

\section{Left invariant Hamiltonian systems and the Maximum Principle}

 Consider now the necessary
conditions of optimality provided by the Maximum Principle. \ The Maximum
Principle states that each minimizer is the projection of an extremal curve
in the cotangent bundle $T^{\ast }G\ $ and each extremal curve is
an integral curve of a certain Hamiltonian vector field on $T^{\ast }G.\ $\
To state all this  in more detail requires \ additional notation and terminology.

As already stated earlier, $\mathfrak{g}^{\ast }$\ denotes the dual of $%
\mathfrak{g.}$ The dual of a Lie algebra\ carries a Poisson structure  inherited from the symplectic structure of $T^{\ast
}G$ realized as the product  $G \times  \mathfrak{g}^{\ast }$ via the left translations. \ \ Functions on $\mathfrak{g}^{\ast }$%
\ are called left-invariant\ Hamiltonians.\ If $f$\ and $h\ $are
left-invariant Hamiltonians then their Poisson bracket $\{f,h\}$ is defined
by $\{f,h\}(l)=l([df,dh]),$\ for\ $l\in \mathfrak{g}$.\ 

On  semisimple  Lie algebras $\mathfrak{g}%
^{\ast }$ can be  identified$\ $with\ $\mathfrak{g}$ via
the quadratic form $\langle\ ,\ \rangle$\ with \ $\langle L , X\rangle =l(X)\ $for all $X\in $\ $%
\mathfrak{g}. $ In this identification  $\mathfrak{p}^{\ast }$\ and $\mathfrak{k}^{\ast }\ $  are  identified with $\fp$ and $\fk$ whenever $\fg$ admits a Cartan decomposition $\fg=\fp\oplus\fk$. The above then implies that $l=l_{\mathfrak{p}}+l_{%
\mathfrak{k}}\ $with\ $l_{\mathfrak{p}}\in \mathfrak{p}^{\ast }\ $and\ $l_{%
\mathfrak{k}}\in \mathfrak{k}^{\ast }$ $\ $is identified with$\mathfrak{\ }%
L=L_{\mathfrak{p}}+L_{\mathfrak{k}}$\  where  $l_{%
\mathfrak{p}}\ $and\ $l_{\mathfrak{k}}$ correspond to $L_{%
\mathfrak{p}}\in $ $\mathfrak{p}$\ and\ $L_{\mathfrak{k}}\in $ $\mathfrak{k.}
$

To preserve the left
invariant symmetries $ T^{\ast }G $ will be  trivialized by the left
translations and considered as the product $G\times $ $\mathfrak{g}^{\ast }.$
 The advantage of the above choice of trivialization is that the Hamiltonian lift of a left invariant vector field becomes a linear function 
on $\mathfrak{g}^{\ast }.$ Recall that a Hamiltonian lift
 of a vector field $X\ $on a manifold $M\ $%
is a function $H_{X}\ $on $T^{\ast }M\ $defined as $H_{X}(\xi )=\xi (X(x))$\
for each $\xi \in T_{x}^{\ast }M\ $(\cite{jurd}). If $X(g)=gA$ is a left invariant vector field on a Lie group $G$, then $H_X(g,l))=l(A),\, l\in\fg^*$.

 Any left invariant function $h\ 
$generates a Hamiltonian vector field\ $\vec{h}\ $on $G\times \mathfrak{g}%
^{\ast }\ $whose integral curves  $(g(t),l(t))\ $ are the solutions of the
following differential equations%
\begin{equation}
\frac{dg}{dt}=g(t)dh(l(t)),\frac{dl}{dt}=-ad^{\ast }(dh(l))(l(t)),
\label{linvham}
\end{equation}

where $dh$\ denotes\ the\ differential\ of\ $h\ $considered as a element of $%
\mathfrak{g}$\ under the \ isomorphism $(\mathfrak{g}^{\ast })^{\ast
}\longleftrightarrow \mathfrak{g,}$\ and where$\ ad^{\ast }(L):\mathfrak{g}%
^{\ast }\rightarrow \mathfrak{g}^{\ast }$\ is defined by \ $ad^{\ast
}(L)(l)=l\circ ad(L).\ $

 With these notations at our disposal we now apply the Maximum Principle to the affine problem.  The affine problem defines  "cost-extended system" in $\IR\times G$: 
\begin{equation}
\frac{dx}{dt}=\frac{1}{2}||U(t)||^{2},\ \frac{dg}{dt}(t)=g(t)(A+U(t)),\,U(t)\in\fk.
\end{equation}%
The Hamiltonian lift of the cost-extended system is given by:
\begin{equation}
H_{U}(\lambda ,l)=\lambda \frac{1}{2}||U||^{2}+l(A+U),\lambda\in\IR,\,l\in\fg^*.
\end{equation}
The above is a function on $T^{\ast }(\mathbb{R\times }G)$ tivialized as ($%
\mathbb{R\times R)\times (}G\times $ $\mathfrak{g}^{\ast })$ with coordinates $(x,\lambda,g,l)$.$\ $
Each \ control function $U(t)\ $generates a time varying Hamiltonian $%
H_{U(t)}(\lambda ,l)\ $; the
integral curves $\xi (t)\ =(x(t),\lambda (t),g(t),l(t))\ $of the associated Hamiltonian vector field $\vec{H}%
_{U(t)}\ $are the solutions of 
\begin{equation}\label{Hamlft}
\frac{dx}{dt}=\frac{\partial H_{U(t)}}{\partial \lambda },\frac{d\lambda }{dt%
}=-\frac{\partial H_{U(t)}}{\partial x},\frac{dg}{dt}=g(A+U(t)),\frac{dl}{dt}%
=-ad^{\ast }(A+U(t))(l(t)
\end{equation}%
It follows that $\lambda \ $is constant for any solution\ $\xi (t)\ $ since $%
\frac{\partial H_{U(t)}}{\partial x}=0.$

\begin{proposition}
$\ $\textbf{The Maximum Principle.} \ Assume that$\ \bar{U}(t)\ $ is an
optimal control that generates the trajectory $\bar{g}(t)$. Let $\bar{x}(t) $ denote its running cost $\ \int_{0}^{t}\frac{1}{2}||U||^{2}dt.\ $\
Then $(\bar{x}(t),\bar{g}(t))\ $ is the projection of an integral curve of $
\bar{\xi}(t)\ =(\bar{x}(t),\bar{\lambda},\bar{g}(t),\bar{l}(t))$\ of $\vec{H}
_{\bar{U}(t)}\ $that satisfies the following conditions:
\begin{equation}
\bar{\lambda}\leq 0. \mbox{ When }\lambda=0  \label{cond1}\mbox{ then, }\bar{l}(t)\neq 0.
\end{equation}
\begin{equation}
\vec{H}_{\bar{U}(t)}(\bar{\lambda},\bar{l}(t))\geq H_{U}(\bar{\lambda},\bar{l
}(t))),U\in \mathfrak{k},a.e. \ in\ [0,T].  \label{cond2}
\end{equation}
\end{proposition}

In the literature on optimal control it is customary to consider only the projections $(g(t),l(t))$ of  integral curves $\xi (t) $ of $\vec{H}_{U}(\lambda ,l)\ \ $  which are parametrized by a non-positive parameter $\lambda$. Control functions $U(t)\ $on\ $T^{\ast }G\ $are called \emph{%
extremal\ }if they generate solutions of (\ref{Hamlft}) that satisfy conditions (\ref{cond1}) and (\ref{cond2}) of the
Maximum Principle.\ \ Extremal curves that correspond to $\lambda =0\ $are
called \ $\emph{abnormal}$ and those that correspond to $\lambda <0\ $are
called \emph{normal.\ \ }In the normal case $\lambda \ $is \ reduced to $-1\ 
$because of the homogeneity properties of $H_{U}(\lambda ,l)$ with respect to $\lambda$. 

The Maximum
Principle can be restated  in terms of the extremals  by saying that each optimal trajectory is the
projection of an extremal curve (normal or abnormal).
Thus the Maximum Principle \ identifies two distinct Hamiltonians
associated with each optimal control problem depending \ on the value of $%
\lambda .$ 

Return now to the Hamiltonians of the affine problem.
\ After the identifications of $l=l_{\mathfrak{p}}+l_{\mathfrak{k}}\ $with\ $%
\mathfrak{\ }L=L_{\mathfrak{p}}+L_{\mathfrak{k}}\ $the Hamiltonian$\
H_{U}(\lambda ,l)\ $\ is identified with $H(\lambda ,L)=\lambda \frac{1}{2}%
||U||^{2}+\langle A,L_{\mathfrak{p}}\rangle+\langle U,L_{\mathfrak{k}}\rangle.$

In the normal case\ \ $\lambda =-1,\ $and the maximality condition (\ref%
{cond2}) easily implies that each \ normal extremal curve\ $(g(t),L(t))\ $
is an integral curve of the Hamiltonian \ \ 
\begin{equation}
H=\frac{1}{2}||L_{\mathfrak{k}}||^{2}+\langle A,L_{\mathfrak{p}}\rangle  \label{NorHam}
\end{equation}%
generated by the extremal control $U(t)=L_{\mathfrak{k}}(t).\ $This
Hamiltonian will be referred to as \emph{the  affine Hamiltonian.}

In the abnormal case \ the maximization relative to $U$ \ \ results in a
constraint $L_{\mathfrak{k}}=0$\ and does not directly yield the value for $U.\ $ Further investigations of these extremals  will be deferred to the next section.
\subsection{Extremal equations}

\ 

 Hamiltonian equations (4.1)  reveal their symmetries more readily when recast on  the Lie algebras rather than on their  duals. In order to  treat the affine Hamiltonian both as  a Hamiltonian on $G$ and as a Hamiltonian on the semidirect product $G_s$, equations (4.1) need to be recast  on  $\fg$ and $\fg_s$. Since the Lie bracket is
different in two cases the differential equations take on different forms.

Recall that $[A,B]_{s }$ denotes the Lie bracket that deforms the
semisimple Lie bracket when $s =1$ to the semidirect  Lie bracket when $%
s =0.$\ Let $dh=dh_{\mathfrak{p}}+dh_{\mathfrak{k}},L=L_{%
\mathfrak{p}}+L_{\mathfrak{k}},X=X_{\mathfrak{p}}+X_{\mathfrak{k}}\ $denote
the appropriate decompositions relative to $\mathfrak{p}$\ and $\mathfrak{k.}
$

Then$\ \frac{dl}{dt}(X)=-ad^{\ast }(dh(l)(l(t)(X)=-l([dh,X]_{s})\ $%
corresponds to $\langle \frac{dL}{dt},X\rangle=-\langle L,[dh,X]_{s }\rangle$
;\ the\ latter\ implies\ that\ \ 

\ $\langle \frac{dL_{\mathfrak{p}}}{dt},X_{\mathfrak{p}}\rangle+\langle \frac{dL_{\mathfrak{k}}%
}{dt},X_{\mathfrak{k}}\rangle=-\langle L_{\mathfrak{p}},[dh_{\mathfrak{p}},X_{\mathfrak{k}%
}]+[dh_{\mathfrak{k}},X_{\mathfrak{p}}]\rangle-\langle L_{\mathfrak{k}},[dh_{\mathfrak{k}%
},X_{\mathfrak{k}}]+s \lbrack dh_{\mathfrak{p}},X_{\mathfrak{p}}]\rangle\
$, or

\ $\langle \frac{dL_{\mathfrak{p}}}{dt},X_{\mathfrak{p}}\rangle+\ \langle\frac{dL_{\mathfrak{k}}%
}{dt},X_{\mathfrak{k}}\rangle=\langle [dh_{\mathfrak{k}},L_{\mathfrak{k}}]+[dh_{\mathfrak{%
p}},L_{\mathfrak{p}}],X_{\mathfrak{k}}\rangle+\langle [dh_{\mathfrak{k}},L_{\mathfrak{p}%
}]+s \lbrack dh_{\mathfrak{p}},L_{\mathfrak{k}}],X_{\mathfrak{p}}\rangle.$

Therefore,\ 
\begin{equation}
\frac{dL_{\mathfrak{k}}}{dt}=[dh_{\mathfrak{k}},L_{\mathfrak{k}}]+[dh_{%
\mathfrak{p}},L_{\mathfrak{p}}],\frac{dL_{\mathfrak{p}}}{dt}=[dh_{\mathfrak{k%
}},L_{\mathfrak{p}}]+s\lbrack dh_{\mathfrak{p}},L_{\mathfrak{k}}].
\label{genham}
\end{equation}

In the case of  affine Hamiltonian $H$ given by (\ref{NorHam}),  $dH=L_{\mathfrak{k}}+A\ $and the preceding
equations become
\ 
\begin{equation}
\frac{dL_{\mathfrak{p}}}{dt}(t)=[L_{\mathfrak{k}},L_{\mathfrak{p}}]+s
\lbrack A,L_{\mathfrak{k}}]=[s A-L_{\mathfrak{p}},L_{\mathfrak{k}}],%
\frac{dL_{\mathfrak{k}}}{dt}(t)\ =[A,L_{\mathfrak{p}}].  \label{norext}
\end{equation}

Abnonormal extremals are the integral curves of\ $H_{0}=\langle L,A+U\rangle \ $subject to the constraint
$L_{\mathfrak{k}}$ $=0,$  that is, abnormal extremal curves $%
(g(t),L(t))$\ are the \ solutions of\ 
\begin{eqnarray}
\frac{dg}{dt} &=&g(t)(A+U(t)),  \label{abn} \\
\ \frac{dL_{\mathfrak{p}}}{dt}(t) &=&[U(t),L_{\mathfrak{p}}]+s
\lbrack A,L_{\mathfrak{k}}],\frac{dL_{\mathfrak{k}}}{dt}(t)\ =[A,L_{%
\mathfrak{p}}]+[U(t),L_{\mathfrak{k}}]
\end{eqnarray}

subject to$\ L_{\mathfrak{k}}(t)=0$\ for all $t\in \lbrack 0,T].$\ 
They are described by the following proposition

\begin{proposition}
Abnormal extremal curves are the solutions of $\frac{dg}{dt}=g(t)(A+U(t))$ generated by
 bounded and measurable controls $U(t) \in \mathfrak{k\ }$\
that satisfy the constraints \ 
\begin{equation}
\lbrack A,L_{\mathfrak{p}}]=0,\ [L_{\mathfrak{p}},U(t)]=0
\end{equation}%
for some element $L_{\mathfrak{p}}\ $in\ $\mathfrak{p}.$

 If $L_\fp$ is regular, then the corresponding abnormal extremal curve whose projection on $G$ is optimal is also normal.\end{proposition}

\begin{proof}
Suppose that $g(t),\ L_{\mathfrak{p}}(t),\ L_{\mathfrak{k}}(t)=0\ $is an
abnormal extremal curve generated by the control $U(t)$\ in $\mathfrak{k}.\ $%
Then equations (\ref{abn}) imply that $\frac{dL_{\mathfrak{p}}}{dt}=[U(t),L_{%
\mathfrak{p}}]\ $and $[A,L_{\mathfrak{p}}]=0.\ $This means that $L_{%
\mathfrak{p}}(t)\ $belongs to the maximal abelian subalgebra $\mathbb{A\ }$%
in $\mathfrak{p}$\ that contains $A.$\ Therefore,\ $\frac{dL_{\mathfrak{p}}}{%
dt}\ $\ also belongs to $\mathbb{A}$.\ 

If $X$ is an arbitrary element of $\mathbb{A\ }$then $\langle [U(t),L_{%
\mathfrak{p}}]\,,X\ \rangle=\langle U(t),[L_{\mathfrak{p}},X]\rangle=0.$\ Therefore, $\langle X,%
\frac{dL_{\mathfrak{p}}}{dt}\rangle=0.\ $Since \ the Killing form is nondegenerate
on $\mathbb{A}$,\ $\frac{dL_{\mathfrak{p}}}{dt}=0$ \ and therefore,$\ L_{%
\mathfrak{p}}$ $(t)\ $\ is constant.\ This proves the first part of the
proposition.

To prove the second part assume that $L_{\mathfrak{p}}\ $is regular. Then\ $%
[A,L_{\mathfrak{p}}]=0\ $\ implies that $[L_{\mathfrak{p}},[A,U(t)]]=0.\ $%
Since $L_{\mathfrak{p}}$\ is regular and belongs to $\mathbb{A},[A,U(t)]\ $%
also belongs to $\mathbb{A}.$ \ It then follows\ that $[A,U(t)]=0\ $by the
argument identical \ to the one \ used in the preceding paragraph.\ 

It now follows that $g(t)=g(0)e^{At}h(t)$ where $h(t)\ $denotes the solution\ of$\ \frac{dh}{dt}(t)=h(t)U(t),\ h(0)=I$.  Let $g_0$ and  $g_1$ denote the boundary points relative to which $g(t)$ is optimal. Then $h(t)$  is optimal relative to $h(0)=I$ and  $h(T)=e^{-AT}g_1$. This means that $h(t)$ is a geodesic in $K$ relative to the bi-invariant metric induced by $\langle\,,\,\rangle$.  Hence the control that generates $h(t)$ must be constant, i.e., $h(t)=e^{Ut}$ for some element $U$ in $\fk$.

The reader can readily verify that  each trajectory $(g(t),U(t))$ of the affine system in which the control $U(t)$ is constant and commutes with $A$ is the projection of a solution of (\ref{norext}).

\end{proof}

\begin{corollary} 
If $\fp$ is such that each non-zero element is regular, then each abnormal extremal that projects onto an optimal trajectory is also a projection of a normal extremal curve. In particular, on isometry groups of space forms (simply connected  symmetric spaces of constant curvature) each optimal trajectory of the affine system is the projection of a normal extremal curve.
\begin{proof}  See the  discussion on space forms  in  Section 6.
\end{proof} 
\end{corollary}
\begin{remark}  The above proposition raises an interesting question. \\
Is every optimal trajectory on an arbitrary symmetric space the projection of a normal extremal curve ?\\
It seems that $G=SL_n(R)$ with $\fp$ the space of symmetric matrices with trace zero and $\fk=so_n(R)$ is a good testing ground for this question. In this situation  there are plenty abnormal extremal curves but it is not clear exactly how they relate to optimality.
\end{remark}

\section{Spectral representation and its consequences}

We will now recall an observation made in (\cite{Rol}\ ) that a system of
differential equations of the form%
\begin{equation}
\frac{dX_{0}}{dt}=[X_{0},X_{1}],\frac{dX_{1}}{dt}=[X_{0},X_{2}],\ldots ,%
\frac{dX_{n}}{dt}=[X_{0},X_{n+1}],\ \frac{dX_{n+1}}{dt}=0  \label{dil}
\end{equation}%
\ admits a spectral representation\ 
\begin{equation}
\frac{dL_{\lambda }}{dt}=[M_{\lambda },L_{\lambda }]  \label{sprep}
\end{equation}

with \ $M_{\lambda }=\frac{1}{\lambda }X_{0}\ $and\ $L_{\lambda }=\frac{1}{%
\lambda }X_{0}+X_{1}+\lambda X_{2}+\cdots \lambda ^{n}X_{n+1}.\ $This
representation is a consequence of a dilational symmetry\ 
\begin{equation}
\tilde{X}_{0}=\frac{1}{\lambda }X_{0},\tilde{X}_{1}=X_{1},\tilde{X}%
_{2}=\lambda X_{2},\ldots ,\ \tilde{X}_{n+1}=\lambda ^{n}X_{n+1}.
\end{equation}

For then $\tilde{X}_{0},\tilde{X}_{1},\ldots ,\tilde{X}_{n+1}\ $also satisfy
( \ref{dil}\ )\ and therefore,\ $\ \frac{d\tilde{L}}{dt}=[\tilde{X}_{0},%
\tilde{L}]\ $where $\tilde{L}=\tilde{X}_{0}+\tilde{X}_{1}+\cdots \tilde{X}%
_{n+1}.\ $

Extremal\ equations\ (\ref{norext}) are of the form (\ref{dil}\ )
with $X_{0}=L_{\mathfrak{p}}-\epsilon A,\ X_{1}=-L_{\mathfrak{k}}\ $and $%
X_{2}=A.\ $Since ( \ref{sprep}) is invariant under a multiplication by $%
\lambda \ $\ it will be convenient to redefine \ $L_{\lambda }\ $as$\
L_{\lambda }=X_{0}+\lambda X_{1}+\lambda ^{2}X_{2}+\cdots \lambda
^{n+1}X_{n+1}\ $in\ which case%
\begin{equation}
L_{\lambda }=L_{\mathfrak{p}}-\lambda L_{\mathfrak{k}}+(\lambda
^{2}-s )A  \label{Spectmat}
\end{equation}%
\ 

is a spectral matrix for equations\ \ (\ref{norext}), in the sense that
 the spectral invariants of $L_{\lambda }\ $are constants of
motion for the corresponding Hamiltonian system. \ Moreover,
 these functions are in involution according to the following proposition.

\begin{proposition}
\label{isospectr}The spectral invariants of $\ L_{\lambda }=L_{\mathfrak{p}%
}-\lambda L_{\mathfrak{k}}+(\lambda ^{2}-1)A\ $Poisson commute with each
other relative to the semisimple Lie algebra structure, while the spectral
invariants of \ $L_{\lambda }=L_{\mathfrak{p}}-\lambda L_{\mathfrak{k}%
}+\lambda ^{2}A\ $Poisson commute relative to the semidirect product
structure.
\end{proposition}

\ \ The proof below is a minor adaptation of the one presented in \ (\cite%
{per}) and (\cite{rey}).\ 

\begin{proof}
Let $T:\mathfrak{g}^{\ast }\rightarrow \mathfrak{g}^{\ast }\ $be defined by $%
T(p+k)=\frac{1}{\lambda }p-k+\mu a\ $for\ $\ p\in \mathfrak{p}^{\ast },k\in 
\mathfrak{k}^{\ast }.$\ Here,\ $a\ $is a fixed element of $\mathfrak{p}%
^{\ast }\ $and \ $\lambda $\ and\ $\mu \ $are parameters. Then, \ $%
T^{-1}=\lambda p-k-\lambda \mu a.\ $This diffeomorphism extends to a
diffeomorphism on forms according to the following formula:
  \begin{equation}
  \{f,g\}_{\lambda ,\mu }(\xi )=(T\circ
\{f,g\})(\xi )=\{f\circ T^{-1}, g\circ T^{-1}\}(T(\xi )),
\end{equation}
 where\ $\{\ ,\
\}\ $denote the canonical Poisson form on $\mathfrak{g}^{\ast }$\ (relative
to the semisimple structure).\ A simple calculation\ shows that 
\begin{equation}
\{f,g\}_{\lambda ,\mu }\ =-\lambda ^{2}\{f,g\}-\lambda \mu
\{f,g\}_{a}-(1-\lambda ^{2})\{f,g\}_{s}  \label{shiftedfms}
\end{equation}

where\ $\{f,g\}_{a}=\{f,g\}(a)$\ and\ $\{f,g\}_{s}\ $is the Poisson bracket
relative to the semidirect product structure. Relative to the
semidirect structure $\{f,g\}_s$ the shifted Poisson bracket $(\{f,g\}_{s})_{\lambda ,\mu
}(\xi )=(T\circ \{f,g\}_{s})(\xi )$ \ takes on a slightly different form:$%
(\{f,g\}_{s})_{\lambda ,\mu }\ =-\{f,g\}_{s}-\lambda \mu (\{f,g\}_{s})_{\ a}$\ 

\ Functions on $\mathfrak{g}^{\ast }\ $which Poisson commute with any other
function on $\mathfrak{g}^{\ast }\ $are called \textit{Casimirs.}, i.e.,
Casimirs are the elements of the center of the Poisson algebra $C^{\infty }(%
\mathfrak{g}^{\ast }).\ $

It \ $f$\ is any Casimir then$\ \ f_{\lambda ,\mu }=T\circ f\ \ $satisfies $%
\{f_{\lambda ,\mu },g\}_{\lambda ,\mu }=0\ $for any function\ $g\ $on \ $%
\mathfrak{g}^{\ast }$and any parameters $\lambda \ $\ and$\ \mu .\ $In the
case that $g$\ is another Casimir then$\ \ \ f_{\lambda ,\mu }\ $and$\
g_{\lambda ,\mu }=T\circ g\ $satisfy$\ \ \ $\ \ 
\begin{equation}
\{f_{\lambda _{1},\mu _{1}},g_{_{\lambda _{2},\mu _{2}}}\}_{\lambda _{1},\mu
_{1}}=\{f_{\lambda _{1},\mu _{1}},g_{_{\lambda _{2},\mu _{2}}}\}_{\lambda
_{2},\mu _{2}}=0.
\end{equation}

for any values $\lambda _{1},\mu _{1},$ $\lambda _{2},\mu _{2}.$The same
applies to the semidirect Poisson bracket.

Suppose now that $\mu \ =\frac{\lambda ^{2}-1}{\lambda }.\ $It follows from
\ (\ref{shiftedfms})\ that\ $\{f,g\}_{\lambda ,\mu }\ =-\lambda
^{2}\{f,g\}-(1-\lambda ^{2})(\{f,g\}_{s}+\{f,g\}_{a}).\ $Therefore,

$\ 0=\frac{1}{\lambda _{1}^{2}-1}\{f_{\lambda _{1},\mu _{1}},g_{_{\lambda
_{2},\mu _{2}}}\}_{\lambda _{1},\mu _{1}}-\frac{1}{\lambda _{2}^{2}-1}%
\{f_{\lambda _{1},\mu _{1}},g_{_{\lambda _{2},\mu _{2}}}\}_{\lambda _{2},\mu
_{2}}$\ $=\frac{\lambda _{2}^{2}-\lambda _{1}^{2}}{(1-\lambda
_{1}^{2})(1-\lambda _{2}^{2})}\{f_{\lambda _{1},\mu _{1}},g_{_{\lambda
_{2},\mu _{2}}}\}.\ $

Since $\lambda _{1},$\ and$\ \lambda _{1}$\ are arbitrary\ $\{f_{\lambda
_{1},\mu _{1}},g_{_{\lambda _{2},\mu _{2}}}\}=0.$ This argument proves the
first part of the proposition because$\ \lambda T=L_{\mathfrak{p}}-\lambda
L_{\mathfrak{k}}+(\lambda ^{2}-1)A\ \ \ $when $\mu \ =\frac{\lambda ^{2}-1}{%
\lambda }\ $after the identifications\ $p\rightarrow L_{\mathfrak{p}},\ $\ \ 
$k\rightarrow L_{\mathfrak{k}}\ $and \ $a\rightarrow A.$

 In the semidirect product \ $\lambda T=L_{\mathfrak{%
p}}-\lambda L_{\mathfrak{k}}+\lambda ^{2}A\ $ when $\lambda =\mu .\ $ Then,

$0=\frac{1}{\lambda _{1}^{2}(}\{f_{\lambda _{1},\mu _{1}},g_{_{\lambda
_{2},\mu _{2}}}\}_{s})_{\ \lambda _{1},\mu _{1}}-\frac{1}{\lambda _{2}^{2}}%
(\{f_{\lambda _{1},\mu _{1}},g_{_{\lambda _{2},\mu _{2}}}\}_{s})_{\ \lambda
_{2},\mu _{2}}=(\frac{1}{\lambda _{2}^{2}}-\frac{1}{\lambda _{1}^{2}}%
)\{f_{\lambda _{1},\mu _{1}},g_{_{\lambda _{2},\mu _{2}}}\}_{s},$

Therefore,\ $\{f_{\lambda _{1},\mu _{1}},g_{_{\lambda _{2},\mu
_{2}}}\}_{s}=0.$
\end{proof}

\section{Specific cases}

It will be convenient to single out  some symmetric pairs $(G,K)$ on which more detailed integrability investigations can be carried out.

\subsubsection{\textbf{Non compact Riemannian symmetric spaces }}

\ \textbf{\ }Every semidirect product $G=V\rtimes H$ admits an involutive
automorphism $\sigma (v,h)=(-v,h),(v,h)\in \ V\ltimes H\ $with $%
K=\{0\}\times H\ $\ the group of fixed points under $\sigma .\ $The\
corresponding splitting is given by $\mathfrak{p}=V\times \{0\}$\ and\ $%
\mathfrak{k}=\{0\}\times \mathfrak{h}.\ $The pair $(G,K)\ $is a Euclidean
symmetric pair. Below are \ some examples of non-Euclidean symmetric spaces.

$\ $\textbf{Selfadjoint groups.\ }A matrix group $G$\ is called \emph{self
adjoint} if the transpose $g^{T}$\ \ belongs to $G$\ \ for every $g\in G.$\
Let $G\ \subseteq SL_{n}(R)\ $be any self adjoint group.\ Define $\sigma
:G\rightarrow G$\ \ by $\sigma (g)=(g^{T})^{-1}.$\ \ Then$\ \sigma (g)=g\ $\
if and only if $g\in SO_{n}(R)\cap G$. Assuming that $G\neq $\ \ $SO_{n}(R)$%
, then $\sigma \ $is an involutive automorphism\ and $(G,K)$\ is a
symmetric pair with $K=G\cap SO_{n}(R).$\ \ The splitting of $\ \mathfrak{g\ 
}$induced\ by$\ \sigma \ $is given by \ 
\begin{equation}
\mathfrak{p}=\{A\in \mathfrak{g}:A^{T}=A\},and\ \mathfrak{k}=\{A\in 
\mathfrak{g}:A^{T}=-A\}.
\end{equation}

The quadratic form defined by 
\begin{equation}
\langle A,B\rangle=\frac{1}{2}Trace(AB)
\end{equation}%
\ is \ $Ad_{K}$\ invariant and positive definite on $\fp$ and negative definite on $\fk$. Hence, $(G,K)\ $is\ a\
Riemannian\ symmetric\ pair. Below are some noteworthy special cases \ of
self adjoint \ groups.

\textbf{Positive definite matrices}.$\ G=SL_{n}(R),\ K=SO_{n}(R).\ $Then$\
SL_{n}(R)/SO_{n}(R)\ $can be identified with\ the\ space\ of\ \ positive
definite $n\times n\ \ $matrices\ with real entries.

\textbf{The generalized upper half plane.\ }\ $G$ $=\ Sp_{n},K=SO_{2n}\cap
Sp_{n}=SU_{n}.$ 

Recall that \ $Sp_{n}\ $denotes$\ $the group that leaves the
symplectic form\ $<x,y>=\sum_{i=1}^{n}x_{i}y_{i+n}-y_{i}x_{i+n}\ $\ in $%
\mathbb{R}^{2n}\ $invariant.\ \ In this situation

\begin{equation*}
\mathfrak{p}=\{\left( 
\begin{array}{cc}
A & B \\ 
B & -A^{T}%
\end{array}%
\right) , B^{T}=B\},\mathfrak{k}=\{\left( 
\begin{array}{cc}
A & B \\ 
-B & A%
\end{array}%
\right) :A^{T}=-A,\ B^{T}=B\}.
\end{equation*}%

The quotient\ $Sp_{n}/SU_{n\ }$can be considered as the generalized upper
half plane since it$\ $can can be realized also as the space of complex $%
n\times n$\ matrices $Z$\ \ with 
\begin{equation}
Z=X+iY
\end{equation}

with \ $X\ and\ Y\ $real $n\times n$\ symmetric matrices and $Y\ $positive
definite.\ 

 For $n=1,$\ $Sp_{1}=SL_{2},\ SU_{1}=SO_{2}(R)\mathcal{\ }$and$%
\ \ \ Sp_{1}/SU_{1\ }\ $coincides with $SL_{2}(R)/SO_{2}(R).$ The latter, with its Riemannian metric induced by $\langle\,,\,\rangle$, is identified with Poincar\'e's upper half plane.

\textbf{Open sets of Grassmannians}$.\ G=$ $SO(p,q),\ K=$ $(O_{p}(R)\times
O_{q}(R))\cap SO_{p+q}(R).$

The quotient space  ${M}_{p,q}^{+}=SO(p,q)/(O_{p}(R)%
\times O_{q}(R))\cap SO_{p+q}(R)\ $\ is identified with the open subset of
 Grassmannians $Gr_{q}(\mathbb{R}^{p+q})\ $consisting of all $q\ $%
dimensional subspaces in \ \ $\mathbb{R}^{p+q}$\ on which the quadratic form
\ \ $\langle x,y\rangle_{p,q}=-\sum_{i=1}^{p}x_{i}y_{i}+\sum_{i=p+1}^{p+q}x_{i}y_{i}$ \
is positive definite. \ The corresponding Lie algebra \ splitting is given by%
\begin{eqnarray*}
\mathfrak{p} &=&\{M=\left( 
\begin{array}{cc}
0 & X^{T} \\ 
X & 0%
\end{array}%
\right) ,\ X\ any\ p\times q\ matrix\},and \\
\mathfrak{k} &=&\{M=\left( 
\begin{array}{cc}
B & 0 \\ 
0 & C%
\end{array}%
\right) :B^{T}=-B,\ C^{T}=-C\}.
\end{eqnarray*}

The space ${M}_{1,n}^{+}\ $\ can be\ identified\
with\ \ the hyperboloid $\mathbb{H}^{n}=\{x\in \mathbb{R}%
^{n+1}:x_{0}^{2}-(x_{2}^{2}+\cdots +x_{n}^{2})=1,\ x_{0}>0\}\ $via\ the
following\ identification.
Let $P$ denote the orthogonal complement\ \ relative to \ $\langle \
,\ \rangle _{1,n}$
of an   $n$\ dimensional subspace $Q$ in \ $\mathbb{R}^{n+1}$ on which\
\ \ $\langle x,y\rangle_{1,n}=-x_{0}y_{0}+\sum_{i=2}^{n}x_{i}y_{i}$ \ is positive
definite. Since $\langle \ ,\ \rangle _{1,n}$\  is positive on $Q$ and non-degenerate on $\mathbb{R}^{n+1}$, 
 $P$\ is transversal to $Q\ $and\ hence\ is one dimensional.\ \ Let\ $%
p=(p_{0},\ldots ,p_{n})\ $be any non-zero point of $P.$\ Since the form $\langle \
,\ \rangle _{1,n}\ $is indefinite on \ \ $\mathbb{R}^{n+1},$\ \ \ $\langle p,p\rangle_{1,n}<0.$\
If \ $p\ $is normalized so that $p_{0}>0$\ $\ $and$\ \langle p\ ,p\ \rangle =-1$\ $\ $%
then\ $p\in $\ $\mathbb{H}^{n}$\ and$\ Q\ \ $is identified with the tangent
space at $p$.

\subsubsection{\textbf{Compact}\ \textbf{Riemannian symmetric spaces}\ }

\paragraph{ \textbf{The Grassmannians.}}

\textbf{\ }Let $G=SO_{p+q}(R)\ $\ with the automorphism $\sigma (g)=JgJ^{-1}$ where\ $J$\ \ denotes a diagonal matrix with its first $p\ $diagonal
entries equal to $1\ $and the remaining diagonal entries equal to $-1.$ It
follows that $\sigma (g)=g\ \ $ if and only if $\ J=gJg^{T}.\ $An easy
calculation shows that $J=gJg^{T}\ $if and only if $g=\left( 
\begin{array}{cc}
g_{1} & 0 \\ 
0 & g_{2}%
\end{array}%
\right) $ where\ $g_{1}\ $is a $p\times p$\ matrix, \ $g_{2}\ $is a $\
q\times q$\ matrix and $g_1=g_1^T$ and $g_2=g_2^T$. Hence the isotropy group $K$ is equal to $(O_{p}(R)\times O_{q}(R))\cap SO_{p+q}(R)$, which will be denoted  by $S(O_{p}(R)\times O_{q}(R)).$\

The tangent map of $\sigma $\ splits $\mathfrak{g=}so_{p+q}(R)$ into $\mathfrak{p }$, the vector space of matrices $%
P=\left( 
\begin{array}{cc}
0 & B \\ 
-B^{T} & 0%
\end{array}%
\right) \ $ where $B$\ is a $\ p\times q\ $matrix, and $\mathfrak{k\ 
}$ the Lie algebra of $K$\ consisting of matrices $Q=\left( 
\begin{array}{cc}
A & 0 \\ 
0 & D%
\end{array}%
\right) $with $A$\ and $D$ antisymmetric.\ 

The pair $(G,K)$\ is a symmetric
Riemannian pair with the metric on $\mathfrak{p}\ $defined by  the quadratic form $%
\langle P_{1},P_{2}\rangle=-\frac{1}{2}Tr(P_{1}P_{2}).$ The homogeneous space $%
Gr_{p}=SO_{p+q}(R)/S(O_{p}(R)\times O_{q}(R))\ $ is the space of  $p$\
dimensional linear subspaces \ in $\mathbb{R}^{p+q}$ and it is a double cover  of $%
SO_{p+q}(R)/SO_{p}(R)\times SO_{q}(R)\ $. The latter is the space of oriented $p$ dimensional  linear subspaces  in $\mathbb{R}^{p+q}$. 

When $p=1$\ and $q=n$, then  the
set of oriented lines in $\mathbb{R}^{n+1}$\ is identified with the sphere $%
\mathbb{S}^{n}$ and the above gives\ 
\begin{equation*}
\mathbb{S}^{n}=SO_{n+1}(R)/\{1\}\times SO_{n}(R).
\end{equation*}

$\ $\textbf{Complex symmetric matrices.\ }$Let \ G=SU_{n}$  with  $\ \sigma
(g)=(g^{T})^{-1}$.\ \ It follows that $\ \sigma (g)=g$\ if and only if $g\in
SO_{n}(R).$The corresponding splitting of $su_n$ identifies $\mathfrak{k}$ with the real part of matrices in $su_n$ and 
 $\mathfrak{p} $ with the imaginary matrices in $su_n$, i.e., $\fp=\{iY:\ Y=Y^{T}\}$ and $\fk=\{X: X^T=-X\}$.
 
 The pair $%
(SU_{n},SO_{n}(R))\ $is a symmetric Riemannian pair with the metric induced
by the trace form $\langle A,B\rangle=-\frac{1}{2}Tr(AB).$
\ \ The quotient space $M=SU_{n}/SO_{n}(R)\ $can be identified with \
complex matrices of the form $e^{iA}$\ for some symmetric real matrix $A,$  
 because every matrix $\ g\in SU_{n}\ $can be written \ in its polar form as $%
g=e^{iA}R$\ for some $R\in SO_{n}(R).$ 

For \ $n=2,\ M\ \ $is a two
dimensional sphere as can be verified by the following argument. \ If $%
A=i\left( 
\begin{array}{cc}
a & b \\ 
b & -a%
\end{array}%
\right) \ $denote a matrix in $\fp\ $ then $A^{2}=-(a^{2}+b^{2})I$ and
therefore,$\ $%
\begin{equation*}
e^{iA}=I\cos \sqrt{a^{2}+b^{2}}+\frac{i}{\sqrt{a^{2}+b^{2}}}A\sin \sqrt{%
a^{2}+b^{2}},
\end{equation*}%
$\ $ or\ 
\begin{equation*}
e^{iA}=\left( 
\begin{array}{cc}
\cos \sqrt{a^{2}+b^{2}}+\frac{ia}{\sqrt{a^{2}+b^{2}}}\sin \sqrt{a^{2}+b^{2}}
& \frac{ib}{\sqrt{a^{2}+b^{2}}}\sin \sqrt{a^{2}+b^{2}} \\ 
\frac{ib}{\sqrt{a^{2}+b^{2}}}\sin \sqrt{a^{2}+b^{2}} & \cos \sqrt{a^{2}+b^{2}%
}-\frac{ia}{\sqrt{a^{2}+b^{2}}}\sin \sqrt{a^{2}+b^{2}}%
\end{array}%
\right).
\end{equation*}
Then
 $\ x=\cos \sqrt{a^{2}+b^{2}},\ y=\frac{a}{\sqrt{a^{2}+b^{2}}}\sin 
\sqrt{a^{2}+b^{2}}\ and\ z=\frac{b}{\sqrt{a^{2}+b^{2}}}\sin \sqrt{%
a^{2}+b^{2}},$ identifies the above matrix with the sphere $\ x^{2}+y^{2}+z^{2}=1.$ The decomposition\ $g=e^{iA}R\ $%
\ corresponds to the Hopf fibration$~S^{3}\rightarrow S^{2}\rightarrow
S^{1}. $

\subsection{\textbf{Space forms}.}

Simply connected Riemannian spaces of constant curvature, known as \emph{%
space forms,}\ \ consist of hyperboloids (spaces of negative curvature),
spheres ( spaces of negative curvature) and Euclidean spaces (spaces of zero
curvature). The normalized prototypes are the unit hyperboloid $\mathbb{H}%
^{n},\ $the unit sphere $\mathbb{S}^{n}$ $\ $and the Euclidean space $%
\mathbb{E}^{n}.$\ It\ follows\ from\ above\ that%
\begin{equation}
\mathbb{S}^{n}=SO_{n+1}/K,\ \mathbb{H}^{n}=SO(n,1)/K,\ \mathbb{E}^{n}=%
\mathbb{R}^{n}\ltimes SO_{n}(R)/K,  \label{spforms}
\end{equation}

where $\ K=\{1\}\times SO_{n}(R).$

The splitting of the corresponding algebras can be described \ in terms of
the curvature parameter $\epsilon=\pm 1,0\ $with%
\begin{equation}
\mathfrak{p}_{\epsilon}=\{\left( 
\begin{array}{cc}
0 & -\epsilon  p^{T} \\ 
p & 0%
\end{array}%
\right) ,\ p\in \mathbb{R}^{n}\},\mathfrak{\ k}_{\epsilon}=\{\left( 
\begin{array}{cc}
0 & 0 \\ 
0 & X%
\end{array}%
\right) ,\ X\in so_{n}(R)\}.  \label{cstcurv}
\end{equation}%
$\ $
 It will be convenient to introduce a shorthand notation and write 
 \begin{equation}
 M_\ep=G_\ep/K , K=\{1\}\times SO_n(R),
 \end{equation}
 with $G_\ep$ equal to $SO_{n+1}(R)$ when $\ep=1$, $SO(n,1)$ when $\ep=-1$, and $SE_n$ when $\ep=0$.
\section{Affine problem on space forms}

\subsection{ Elastic curves and the pendulum.\ }

\ \ On space forms equations \ (\ref{norext})\ \ admit
additional integrals of motion  in involution with the spectral ones \
described by Proposition \ref{isospectr}.\ 
 
 They are described as follows: let
$\mathfrak{k}_{A}=\{M\in \mathfrak{k:[}M,A]=0\}$,
and let $\mathfrak{k}_{A}^{\perp }$ denote the orthogonal complement of $%
\mathfrak{k}_{A}$\ in $\mathfrak{k}$\ relative to $\langle\,,\,\rangle$. It\ \
is easy to see that $\mathfrak{k}_{A}$\ is a Lie subalgebra of $\mathfrak{k\ 
}$and that $[A,$ $L_{p}]\in $ $\mathfrak{k}_{A}^{\perp }.\ $Therefore, the
projection of $L_{\mathfrak{k}}$\ on $\mathfrak{k}_{A}\ $is constant along
the solutions of (\ref{norext}).\ \ 

Recall now that an extremal control $U(t)\ $\ is equal to $L_{\mathfrak{k}%
}(t)\ $\ and that the \ corresponding extremal energy is equal \ to $\frac{1%
}{2}\int_{0}^{T}||U(t)||^{2}dt.\ \ $Let $L_{\mathfrak{k}}\mathfrak{(t)=}%
L_{A}+L_{A}^{\bot }(t)$\ denote the decomposition of \ $L_{\mathfrak{k}}%
\mathfrak{(t)\ }$onto the factors $\mathfrak{k}_{A}\ $and $\mathfrak{k}%
_{A}^{\perp }.$\ Then, $\ $%
\begin{equation*}
\frac{1}{2}\int_{0}^{T}||U(t)||^{2}dt=\frac{1}{2}%
\int_{0}^{T}(||L_{A}||^{2}+||L_{A}^{\bot }(t)||)dt=\frac{1}{2}%
\int_{0}^{T}(||L_{A}^{\bot }(t)||^{2}+\ constant\ \ 
\end{equation*}%
along each trajectory of (\ref{norext}).\

 Remarkably, \ $%
||L_{A}^{\bot }(t)||^{2}\ =\kappa(t) ^2$, where $\kappa (t)$ denotes the geodesic curvature of the
projected curve on $G_\ep/K$, whenever  $L_A=0$ and $||A||=1$ (  the norm  of a matrix $A=\left (\begin{matrix} 0&-\ep a^T\\a&0\end{matrix}\right )$  is given  by $\sqrt{\sum_{i=1}^n a_i^2}\,)$. 

To demonstrate this fact, consider $M_\ep$  as a principal $G_\ep$ bundle with connection $\mathcal{D}$ consisting of left invariant vector fields on $G_\ep$ that take values in $\fp_\ep$ at the group identity.
In this setting curves $g(t)$in $G_\ep$ are called horizontal if $\frac{dg}{dt}\in \mathcal{D}(g(t))$ or, equivalently, if $g^{-1}(t)\frac{dg}{dt}(t)\in\fp_\ep$ for all $t$.  Curves $x(t)\ $in\ $G_{\ep}/K$ \ can be represented by horizontal curves via the  formula
\[ 
x(t)=g(t)e_{0}, 
e_{0}=\left( 
\begin{array}{c}
1 \\ 
0 \\ 
\vdots\\0%
\end{array}%
\right) \in \mathbb{R}^{n+1}. 
\]  In this representation,
 $||\frac{dx}{dt}||_\ep$,  the Riemannian length in $M_\ep$ of the tangent vector $\frac{dx}{dt}$,  is given by 
$||A_\ep(t)||,\mbox{ where }A_\ep (t)=\frac{dg}{dt}(t)g^{-1}(t).$

 Solution curves $g(t)$ of the affine system \ $\frac{dg}{dt}%
(t)=g(t)(A+U(t))\ \ $project\ onto the same curve\ $x(t)\ $as the
associated horizontal curves $\tilde{g}(t)=g(t)h(t),$ where $h(t)$ is a \
solution\ in\ $K$ of $\frac{dh}{dt}(t)=h(t)U(t).$\ It follows that\ $\frac{d%
\tilde{g}}{dt}(t)=\tilde{g}(t)(h(t)Ah^{-1}(t))$. Hence,

\begin{equation*}
||\frac{dx}{dt}||_\ep=||h(t)Ah^{-1}(t)||=||A||=1.
\end{equation*}
Then, $\frac{D}{dx}(\frac{dx}{dt})$,\
the covariant derivative of $\frac{dx}{dt}$\ along $x(t),\ $is\ given\ by$\ $%
\ \ \ 
\begin{equation*}
\ \frac{D}{dx}(\frac{dx}{dt})=(\tilde{g}(t)(h(t)[U(t),A]h^{-1}(t))\ e_{0}.
\end{equation*}

Since $||\frac{dx}{dt}||_\ep=||A||=1$, the geodesic curvature $\kappa (t)\ $of$\
x(t)\ $is given by 
\begin{equation*}
\kappa ^{2}(t)=|||\frac{D}{dx}(\frac{dx}{dt})||^{2}=||[U(t),A]||^{2}
\end{equation*}

In particular along the extremal curves $U(t)=L_{\mathfrak{k}}$\ hence,\ \ $%
\kappa ^{2}(t)=||[L_{A}^{\bot }(t),A]||^{2}$ when $L_A=0$.

It remains to show that $||[L_{A}^{\bot }(t),A]||^{2}=||L_{A}^{\bot
}(t)||^{2}.$\ There exists $h\in K\ $such that $h^{-1}Ah=E_{1}=\left( 
\begin{array}{cc}
0 & -ke_{1}^{T} \\ 
e_{1} & 0%
\end{array}%
\right)  $,  a consequence of the fact that$\ \ K\ $\ acts transitively
 by adjoint action  on the unit sphere in $\fp_\ep .\ $Then,\ 
\begin{equation*}
0=h^{-1}[A,\mathfrak{k}_{A}]h=[E_{1},h^{-1}\mathfrak{k}_{A}h],
\end{equation*}%
and therefore,%
\begin{equation*}
\ h^{-1}\mathfrak{k}_{A}h=\left( 
\begin{array}{ccc}
0 & 0 & 0 \\ 
0 & 0 & 0 \\ 
0 & 0 & so_{n-1}(R)%
\end{array}%
\right) ,h^{-1}\mathfrak{k}_{A}^{\bot }h=\{\left( 
\begin{array}{ccc}
0 & 0 & 0 \\ 
0 & 0 & -l^{T} \\ 
0 & l & 0%
\end{array}%
\right) ,\ l\in \mathbb{R}^{n-1}\}.
\end{equation*}

Hence,%
\begin{equation*}
\ \kappa ^{2}(t)=||[L_{A}^{\bot }(t),A]||^{2}=||h^{-1}L_{A}^{\bot
}(t)h,E_{1}||^{2}=||e_{1}\wedge l||^{2}=||l||^2.
\end{equation*}%
Therefore,
$||l||^{2}=||L_{A}^{\bot }(t)||,$\ and the extremal
energy $\frac{1}{2}%
\int_{0}^{T}(||L^\perp_{A}||^{2}dt$
is given by 
\begin{equation*}
\frac{1}{2}\int_{0}^{T}\kappa ^{2}(t)dt.\end{equation*}

 Curves $x(t)\ $in the base space $G_{k}/K\ $which are the
projections of these extremal curves are called \emph{elastic\ }and $\frac{1%
}{2}\int_{0}^{T}\kappa ^{2}(t)dt\ $is called their\emph{\ elastic energy }(%
\cite{JF}). Equations\ (\ref{norext})\ with $L_{A}=0$
can  be\ rephrased as
\begin{equation}
\frac{dL_{\mathfrak{p}}}{dt}(t)=[L_{A}^{\bot },L_{\mathfrak{p}}]+s
\lbrack A,L_{A}^{\bot }],\frac{dL_{A}^{\bot }}{dt}(t)\ =[A,L_{\mathfrak{p}}].
\label{Elastic}
\end{equation}

We will return to these equations after a  brief digression to mechanics and the connections  between the elastic problem and the motions of a mathematical pendulum.

\subsubsection{\textbf{The pendulum}.\ }

There is a remarkable\ (and somewhat mysterious)\ connection\ between
elastic curves and \ heavy tops\ that will be recalled below\ (\ see also\ (%
\cite{Mem}) for a more general discussion). \ Consider first an $n\ $%
dimensional pendulum of unit length\ suspended at the origin of $\mathbb{R}%
^{n}$\ and acted upon by the " gravitational force" $\vec{F}=-e_{1},$where $%
e_{1},\ldots ,e_{n}\ $denote the standard basis in $\mathbb{R}^{n} $(here,
 all physical constants are normalized to one). 

The motions of the pendulum are\ confined to\ the unit sphere $\mathbb{S}%
^{n-1}.\ $For each curve curve $q(t)\ $on $\mathbb{S}^{n-1}$\ let$\ \
f_{1}(t),\ldots ,f_{n}(t)$\ denote an orthonormal frame , called the moving frame,
adapted to $q(t)$ by the constraint $q(t)=f_{1}(t)$ and positively oriented relative to the absolute frame $e_{1},\ldots ,e_{n}$,  
 in the sense that  the matrix $R(t)$ defined by $f_{i}(t)=R(t)e_{i}$,\ 
$i=1,\ldots ,n,\ $ belongs to $ SO_{n}(R).$

This choice of polarization identifies the sphere as the quotient $G/K$\ \
with\ $G=SO_{n}(R)\ $\ and $K\ $\ the\ isotropy\ subgroup\ of $SO_{n}(R)\ $%
defined\ by\ $Ke_{1}=e_{1}.\ $Evidently, $K=\{1\}\times SO_{n-1}(R).\ $Let $%
\mathfrak{k}_{0}$\ denote the Lie algebra of $K\ $and let $\mathfrak{k}_{1}%
\mathfrak{\ }$denote the orthogonal complement in $\mathfrak{g}=so_{n}(R)$\
relative to the trace form.\ Then%
\begin{equation*}
\mathfrak{k}_{1}=\{\left( 
\begin{array}{cc}
0 & -u^{T} \\ 
u & 0%
\end{array}%
\right) ,u\in \mathbb{R}^{n-1}\},\mathfrak{k}_{0}\mathfrak{=}\left( 
\begin{array}{cc}
0 & 0 \\ 
0 & so_{n-1}(R)%
\end{array}%
\right) .
\end{equation*}

We will regard $SO_n(R)$ as the principal $SO_n(R)$ bundle ( under the right action) over $\mathbb{S}^{n-1}$ with a connection  $\mathcal{D}$\  consisting of  the left invariant vector fields with values in $%
\mathfrak{k}_1$.   As usual, vector fields in $\mathcal{D}$ and their integral curves will be  called horizontal. 

 It follows that 
every curve $q(t)\ $on \ $\mathbb{S}^{n-1}\ $can be lifted to a horizontal
curve$\ R(t)\ $, in the sense that $q(t)=R(t)e_{1}$, and $\frac{dR}{dt}%
=R(t)\left( 
\begin{array}{cc}
0 & -u^{T}(t) \\ 
u(t) & 0%
\end{array}%
\right) \ $for some curve $u(t)\ $in\ $\mathbb{R}^{n-1} $.  Furthermore, it follows that any two such
liftings are related by a left multiple by an element in $K.$

The kinetic energy $T\ $associated with a path in $\mathbb{S}^{n-1}$\ is
given by 
\begin{equation*}
T=\frac{1}{2}||\frac{dq}{dt}||^{2}=\frac{1}{2}||\frac{dR}{dt}e_{1}||^{2}=%
\frac{1}{2}||R(t)\left( 
\begin{array}{c}
0 \\ 
u(t)%
\end{array}%
\right) ||^{2}=\frac{1}{2}||u(t)||^{2}.
\end{equation*}%
The potential energy $V(q)\ $ relative to a fixed point  point $q_0$ is given by $V(q)=-\int_{q_{0}}^{q}\vec{F}%
\cdot \frac{d\sigma }{dt}dt,\ $where\ $\sigma (t)$\ is a\ path from $q_{0}$\
to\ $q.$\ It follows that\ $V(q)=e_{1}\cdot (q-q_{0}).\ $It\ is\ convenient\
to\ take\ $q_{0}=-e_{1}\ $in which case $V=e_{1}\cdot q+1.\ $

The Principle of Least Action states that each motion $q(t)$ of the pendulum
minimizes the action\ $\int_{t_{0}}^{t_{1}}\mathcal{L(}q(t),\frac{dq}{dt}%
)dt\ $over the paths from $q(t_{0})\ $to $q(t_{1})\ $for any $t_{0}$\ and $%
t_{1}$\ (sufficiently near each other to avoid conjugate points), where $\mathcal{L}$\ denotes the Lagrangian $%
\mathcal{L}=T-V.$\ Thus motions of the pendulum can be viewed as the solutions of the
following optimal control problem\ on $SO_{n}(R)\ $:\ 

Minimize $\ $the$\ $integral\ $\int_{t_{0}}^{t_{1}}(\frac{1}{2}%
||u(t)||^{2}-(e_{1}\cdot Re_{1}+1))dt\ $over the solutions $R(t)\in SO_{n}(R)\ 
$\ of $\frac{dR}{dt}=R(t)U_{1}(t),\ U_{1}(t)=\left( 
\begin{array}{cc}
0 & -u^{T}(t) \\ 
u(t) & 0%
\end{array}%
\right) \ $subject to the boundary conditions $R(t_{0})\in \{R:\
Re_{0}=q_{0}\},R(t_{1})\in \{R:\ Re_{0}=q_{1}\}.$

The Maximum Principle then leads to the\ energy Hamiltonian $\mathcal{H}\ $%
on the cotangent bundle of $SO_{n}(R).\ $\ In \ the realization of the cotangent bundle as the product $SO_n(R)\times so^*_n(R)$, further identified with the tangent bundle $SO_{n}(R)\times so_{n}(R)\ $ via the trace
form $\langle A,B\rangle=-\frac{1}{2}Tr(AB),$ the$\ $energy Hamiltonian is given by $\mathcal{H}=%
\frac{1}{2}\langle Q_{1},Q_{1}\rangle +(e_{1}\cdot Re_{1}+1)\ $where  $Q_{1}$\  denotes
the projection on $\mathfrak{k}_{1}\ $of a matrix $Q\ $in\ $so_{n}(R).$

\ The Hamiltonian equations associated with $\mathcal{\vec{H}\ (}$see (\cite%
{Mem}, Ch. IV) for details)\ are given by :%
\begin{equation}
\frac{dR}{dt}=R(t)(Q_{1}(t)),\frac{dQ}{dt}(t)=[Q_{1}(t),Q(t)]-R^{T}(t)e_{1}%
\wedge e_{1}  \label{Pend}
\end{equation}

It is evident from (\ \ref{Pend}\ ) that the projection $Q_{0}$\ of $Q(t)\ $%
on $\mathfrak{k}_{0}$\ is constant.\ But then this constant must be zero
because of the transversality condition imposed by the Maximum principle. To  be more explicit,  recall that the transversality condition states that each  extremal curve $(R(t),Q(t))\ $ annihilates\ the$\ $%
tangent vectors\ of $S_{0}\ $at $R(t_{0}).$\ Therefore,\ $\langle Q(t_{0}),X\rangle=0\ \ $%
for all $X\in \mathfrak{k}_{0}\ $(\ since the tangent space at\ $R(t_{0})\ $%
\ is equal to $\{R(t_{0})X:\ X\in \mathfrak{k}_{0}\}$). The transversality condition at the terminal time $t_1$ reaffirms that $Q_0=0$; hence, it  and is redundant in this case.

Equations (\ref{Pend}\ ) can be lifted to the semidirect product $\mathbb{R}%
^{n}\ltimes so_{n}(R)\ $by identifying  vector $p(t)=-R^{T}(t)e_{1}$ with matrix  \ $P(t)=\left( 
\begin{array}{cc}
0 & 0 \\ 
p(t) & 0%
\end{array}%
\right)  $ and vector $e_1$ with matrix $\ E_{1}=\left( 
\begin{array}{cc}
0 & 0 \\ 
e_{1} & 0%
\end{array}%
\right) .$ \ Both matrices belong to the Cartan space\ $\mathfrak{p}$ in the
 semidirect  Lie algebra  $\{\left( 
\begin{array}{cc}
0 & 0 \\ 
x & X%
\end{array}%
\right) :\ x\in \mathbb{R}^{n},\ X\in so_{n}(R)\}. $

Let $\tilde{Q}=\left( 
\begin{array}{cc}
0 & 0 \\ 
0 & Q%
\end{array}%
\right) \ $\ denote the embedding of $so_{n}(R)$\ into the semidirect
product algebra, and let $\tilde{Q}_{0}\ $and $\tilde{Q}_{1}\ $denote the
embeddings of $Q_{0}\ $and $Q_{1}.$

Then 
\begin{equation}
\left( 
\begin{array}{cc}
0 & 0 \\ 
Q_{1}p(t) & 0%
\end{array}%
\right) =[P(t),\tilde{Q}_{1}],\ and\ [E_{1},P(t)]=\left( 
\begin{array}{cc}
0 & 0 \\ 
0 & p(t)\wedge e_{1}%
\end{array}%
\right) .
\end{equation}

It follows that the extremal equations (\ \ref{Pend}\ )\ can be written also
as 
\begin{equation}
\frac{dg}{dt}=g(t)(E_{1}+\tilde{K}_{1}(t)),\frac{dP}{dt}=[P(t),\tilde{Q}%
_{1}],\frac{d\tilde{Q}_{1}}{dt}=[E_{1},P(t)],  \label{Pend1}
\end{equation}

where $g(t)=\left( 
\begin{array}{cc}
1 & 0 \\ 
q(t) & R(t)%
\end{array}%
\right) .\ $

 The Lie algebra  part of equations (\ref{Pend1}) agree with equations (\ \ref{Elastic}\
)$\ $when $s =0$,\ \ $\epsilon=1$\ and $A=E_{1}.\ $ So on the level of Lie algebras, the equations of the mathematical pendulum coincide with the equations for  the Euclidean elastic curves.

Consider now the isospectral matrix $\ L_{\lambda }=L_{\mathfrak{p}}-\lambda
L_{A}^{\perp }+(\lambda ^{2}-\epsilon )A $ associated\ with\ (\ %
\ref{Elastic}\ )$.$\ Since the spectral invariants of $L_{\lambda }$\ are\
invariant\ under conjugations by elements in $K$\ there is no loss in
generality if $A$\ is taken to be $E_{1}.$ Matrices in $\mathfrak{%
p}_{\ep},\ \ep=\pm 1$\ are of the form $\left( 
\begin{array}{cc}
0 & -\ep p^{T} \\ 
p & 0%
\end{array}%
\right) ~$\ and can be written as\ $\bar{p}\wedge _{\ep}e_{0}\ \ $where\ $\bar{%
p}=\left( 
\begin{array}{c}
0 \\ 
p%
\end{array}%
\right) ,p\in \mathbb{R}^{n}\ $and where%
\begin{equation}
(a\wedge _{\ep}b)x=(b,x)_{\ep}a-(a,x)_{\ep}b,\ with\
(x,y)_{\ep}=x_{0}y_{0}+ep \sum_{i=1}^{n}x_{i}y_{i},
\end{equation}

while matrices in $L_{A}^{\perp }$\ can be written as $\bar{x}\wedge e_{1}$\
with $\bar{x}=\left( 
\begin{array}{c}
0 \\ 
0 \\ 
x%
\end{array}%
\right) ,x\in \mathbb{R}^{n-1}.\ $It follows that the range of $L_{\lambda
}~ $is contained in the linear span of $\bar{p},\bar{x},e_{0},e_{1}.$ An
easy calculation \ shows that the \ characteristic polynomial of the
restriction of $L_{\lambda }\ $to this vector space is given by 
\begin{equation}
\xi ^{4}+c_1\xi ^{2}+c_2=0\ 
\end{equation}
with \[c_1=\ep(\lambda^2-s)+(\lambda^2-s)H+s||L^\perp_A||^2+||L_\fp||^2,  c_2=\lambda^2(||L^\perp_A||^2||L_\fp||^2-||[L^\perp_A,L_\fp]||^2.\]

 The coefficients $c_1$ and $c_2$ reveal  $I_2=||L^\perp_A||^2||L_\fp||^2-||[L^\perp_A,L_\fp]||^2$ as a new integral of motion in addition to the Hamiltonian $\mathcal{H}$ and the Casimir $I_1=s||L^\perp_A||^2+||L_\fp||^2$.
This \ integral of
motion \ admits a nice geometric interpretation \ relative to the underlying
elastic curve$\ x(t)$:
\begin{equation*}I_{2}=\kappa ^{2}(t)\tau (t)
\end{equation*} where $\tau (t)$\
denotes\ the\ torsion\ of\ $x(t)$.\ Then$\ \xi (t)=$\ $\kappa ^{2}(t)\ \ $%
is a solution of 
\begin{equation}
(\frac{d\xi }{dt})^{2}=-\xi ^{3}+4(\mathcal{H-}k)\xi ^{2}+4(I_{1}-\mathcal{H}%
^{2})\xi +4I_{2}
\end{equation}

and hence is solvable in terms of elliptic
functions\ \ (\cite{JF},\ \cite{Mem}\ ).

These geometric identifications suggest an integrating procedure in terms of
the Serret-Frenet frames \ along an elastic curve $x(t).\ $ Let $T,N,B\ \ \ $%
denote the Serret-Frenet triad given by the standard formulas $\ $%
\begin{equation}
\frac{dT}{dt}=\kappa (t)N(t),\ \frac{dN}{dt}=-\kappa (t)T(t)+\tau (t)B(t),\ 
\frac{dB}{dt}=-\kappa (t)N(t).  \label{S-F frame}
\end{equation}
 
The tangent vector $T(t)$ can be identified with $T(t)=(h(t)E_{1}h^{-1}(t))\ $ in the horizontal distribution $\mathcal{D}$ where $h(t)\ $is a solution of $\frac{dh}{dt}%
(t)=h(t)(U(t))\ $with $\ U(t)=L_{A}^{\perp }(t)\ $. Then the normal  and the binormal vectors $N(t)$ and $B(t)$ can be easily  obtained from (\ref{S-F frame}). It was shown first in \cite{Gr} and then in \cite{JF} that $\frac{dB}{dt}$ is contained in the linear span of $T(t),N(t),B(t)$. Hence,
 equations \ (\ref{S-F frame})\ carry complete information about elastic
curves. By analogy, the equations of the
mathematical pendulum are also integrable by an identical procedure.

It can be shown that the general case with\ \ $L_{A}\ $an arbitrary
constant is related to an $n\ $dimensional heavy top with \ equal principal moments of inertia, but these details will not be addressed
here.

\begin{remark}
The affine distribution $\mathcal{D}(g)=\{A+U:U\in \mathfrak{k}\}\ $does not
extend to the elastic problems on more general symmetric spaces because
the isotropy group $K$\ does not act transitively on the spheres in the
Cartan space $\mathfrak{p}$. Hence not every curve$\ x(t)\ $in $G/K\ $can be
lifted to\ a horizontal curve\ $g(t)\ $that is a solution of $\frac{dg}{dt}%
=g(t)((h(t)Ah^{-1}(t)),\ $for some curve $h(t)\ $in\ $K.\ $\ This observation raises a question about the geometric
significance of the affine problem for general symmetric spaces.
\end{remark}

\section{Affine problem on coadjoint orbits}

Certain coadjoint orbits \ coincide with the cotangent bundles of quadric
surfaces and the \ restriction of the affine Hamiltonian to these orbits
coincides with the Hamiltonians associated with mechanical systems with
quadratic potential. \ On these orbits the spectral invariants of (\ref%
{Spectmat}\ )\ form Lagrangian submanifolds\ of the orbits, or,  stated
differently, the restrictions of the Hamiltonian to such manifolds become completely
integrable. These findings \ \ provide a natural theoretical
framework for several classically known integrability results\ and at the
same time point to a larger class of \  systems that conform to the same integration procedures. The  text below supports these
claims in complete detail.

Recall that \ the coadjoint orbit \ $\mathcal{O}_{G}(l_{0})\ $of $G\ \ $%
through\ $l_{0}\in \mathfrak{g}^{\ast }$\ is\ defined
\begin{equation*}
\mathcal{O}_{G}(l_{0})\ =\{l\ :l=Ad_{g^{-1}}^{\ast }(l_{0}),g\in G\}, 
\end{equation*}%
where $Ad_{g^{-1}}^{\ast }(l_{0})(X)=l_{0}(Ad_{g^{-1}}(X)),\ X\in \mathfrak{%
g.}$\ Also recall that$\mathfrak{\ g}^{\ast \ }$is\ a\ Poisson\ manifold\
under the Poisson bracket$\ \{f,h\}(l)=l([df,dh]),\ l\in \mathfrak{g}
^{\ast }$, \ and that $\mathfrak{\ g}^{\ast }$ is foliated by coadjoint
orbits of $G\ $each of which is symplectic. \ More precisely,\ the tangent
space of$\ \mathcal{O}_{G}(l_{0})$\ at $l$\ consists of vectors $v=ad^{\ast
}M(l),\ M\in \mathfrak{g}\ $\ and\ the symplectic form $\omega \ $at $l\ $is
given by 
\begin{equation}
\omega _{l}(v_{1},v_{2})=l([M_{1},M_{2}]),\ with\ \ v_{1}=ad^{\ast
}M_{1}(l)\ and\ v_{2}=ad^{\ast }M_{2}(l).  \label{symp}
\end{equation}

$\ $(\cite{Ar}, Appendix 2).

In the semisimple case $\mathcal{O}_{G}(l_{0})\ $is identified with the
adjoint orbit $Ad_{G}(L_{0})\ $via the correspondence \ $\langle L\ ,\ X\rangle=l(X)\ $%
for all $X\in $\ $\mathfrak{g.}$ Consequently, each adjoint orbit in a semisimple Lie algebra is even
dimensional. In this correspondence, tangent vectors at $l\ $are identified
with\ \ matrices\ $v=$ $[L,M]\ $and the symplectic form\ takes on its dual
form $\omega _{L}(v_{1},v_{2})=\langle L,[M_{1},M_{2}]\rangle.$

 When $(G,K)\ $\ is a symmetric pair  than \ $\mathfrak{g^*}$ carries another Poisson structure $\{\,,\,\}_s$  induced by the  semidirect
product  $\mathfrak{p}\rtimes \fk$.  As in the semisimple case the quadratic form $\langle \,,\,\rangle$ can be used to  identify the coadjoint orbits with certain submanifolds of $\fg$. Since the Killing form is not invariant relative to the semidirect Lie bracket, these manifolds  need not coincide with the adjoint orbits. The proposition below describes their structure.
\begin{proposition}
Suppose that $\ l_{0}\in \mathfrak{g}^{\ast },\ g=(X,h)\in G_{s}\ $and $%
l=Ad_{g-1}^{\ast }(l_{0}),\ $and further suppose that$\ l_{0}\longrightarrow
L_{0}=P_{0}+Q_{0},\ and\ l\longrightarrow L=P+Q$\ \ are the correspondences
defined by the Killing form\ with $P_{0}$\ and $P\ $in $\mathfrak{p}$\ and $%
Q_{0}$\ and $Q\ $in $\mathfrak{k.}$Then%
\begin{equation}
P=Ad_{h}(P_{0}),\ and\ Q=[Ad_{h}(P_{0}),X]+Ad_{h}(Q_{0}).  \label{sdorbit}
\end{equation}

\begin{proof}
\ Let $\ Z=U+V\ $be an arbitrary point of $\mathfrak{g\ }$with\ $U\in 
\mathfrak{p\ }$and $V\in \mathfrak{k.}$

Then$\ $

$Ad_{g^{-1}}(Z)=\frac{d}{d\epsilon }(g^{-1}(\epsilon U,e^{\varepsilon
V})g)|_{\varepsilon =0}=\frac{d}{d\epsilon }(-Ad_{h^{-1}}(X)+Ad_{h^{-1}}(%
\epsilon U+e^{\varepsilon V}(X)e^{-\varepsilon V}),h^{-1}e^{\varepsilon V}h)$

$=Ad_{h^{-1}}(U+[X,V])+Ad_{h^{-1}}(V).$\ \ 

Hence,

$\
l(Z)=l_{0}(Ad_{g-1}(Z)=\langle P_{0},Ad_{h^{-1}}(U+[X,V])\rangle+\langle Q_{0},Ad_{h^{-1}}(V)\rangle $

$=\langle Ad_{h}(P_{0}),U\rangle+\langle Ad_{h}(P_{0}),[X,V]\rangle+\langle Ad_{h}(Q_{0}),V\rangle$

$=\langle Ad_{h}(P_{0}),U\rangle+\langle[Ad_{h}(P_{0}),X],V\rangle+\langle Ad_{h}(Q_{0}),V\rangle =\langle P,U\rangle+\langle Q,V\rangle$

\ Therefore,

$P=Ad_{h}(P_{0}),$\ and $Q=[Ad_{h}(P_{0}),X]+Ad_{h}(Q_{0}).$
\end{proof}
\end{proposition}

\ For left invariant Hamiltonians $H\ $ each coadjoint orbit is an integral manifold \ for the Hamiltonian vector
field $\vec{H}$. Moreover,\ the Hamiltonian vector field on a  coadjoint orbit induced by the restriction\ of $%
H$ 
coincides with the restriction of $\vec{H}\ $ to the coadjoint orbit. The
latter fact will be of central importance for the rest of the paper as it
will be shown that the Hamiltonian \ (\ref{NorHam}) \ restricted to the
coadjoint orbits through rank one matrices in $SL_{n+1}(R)$ relate
directly to Kepler's problem, geodesic problem of Jacobi and the mechanical
problem of Newman. The identification  with the affine problem provides natural explanation for their 
integrability.

\subsection{Coadjoint orbits\ on the vector space of matrices of trace zero}

\bigskip

\ The vector space $\ V_{n}\ $of $n\times n\ $matrices of trace zero admits
\ several kinds of \ Lie algebras \ and each of these Lie algebras induces
its own \ Poisson structure on the dual space$\ V_{n}^{\ast }$. The most
common Poisson structure is \ induced by the canonical \ Lie bracket,\ i.e.,
\ in which $V_{n}$\ as a Lie algebra is equal to \ $sl_{n}(R).$ \ The
decomposition of $sl_{n}(R)\ $as the sum of symmetric and skew-symmetric
matrices ( associated with the automorphism $\sigma (g)=(g^{T})^{-1}\ )\ $%
allows $V_{n}\ \ $to \ be considered also as the semidrect Lie algebra \ 
\emph{Sym}$_{n}\rtimes so_{n}(R),\ $where \emph{Sym}$_{n}\ $denotes the
space of symmetric $n\times n\ $matrices of trace zero. There are also automorphisms of
non-Riemannian type\ which induce semidirect products  of their own.  The paragraph below describes
these semidirect products in some detail.

Let$\ n=p+q\ $and let $\sigma :SL_{p+q}(R)\rightarrow SL_{p+q}(R)$ be defined by $\ \sigma (g)=J((g^{T})^{-1})J^{-1}$
where$\ J$\ \ is diagonal matrix with its first $p\ $diagonal entries equal
to $1\ $and the remaining diagonal entries equal to $-1.\ $ $\ $It follows
that $\ \sigma (g)=g$\ \ if and only if $J=gJg^{T}\ $ or equivalently, $\sigma (g)=g$ if
and only if$\ g\in SO(p,q).$ The$\ $tangent\ map\ $\sigma _{\ast }$\ given by
\ $\sigma _{\ast }(A)=-JA^{T}J$\ induces a decomposition $\mathfrak{p}%
\oplus \mathfrak{k}$\ where \ $\mathfrak{p}$\ consists of matrices$\ P$\
such that $P=JP^{T}J\ $which implies that $P=\left( 
\begin{array}{cc}
A & B \\ 
-B^{T} & D%
\end{array}%
\right) \ $with $A$\ and\ $D\ $symmetric and$\ B\ $an arbitrary $p\times q\ 
$matrix.\ The Lie algebra\ $\mathfrak{k}\ $\ is the Lie algebra of $%
K=SO(p,q)\ $ and consists of matrices $Q=\left( 
\begin{array}{cc}
A & B \\ 
B^{T} & D%
\end{array}%
\right) \ $with $A\ $and $D\ $antisymmetric.\ 

The symmetric pair $\ (SL_{p+q}(R),SO(p,q))\ $is \ strictly pseudo-
Riemannian because the subgroup\ of the restrictions of $Ad_{K}\ $to\ $%
\mathfrak{p}$\ is not a compact subgroup of $Gl(\mathfrak{p})$\ (\ the trace
form is indefinite on $\mathfrak{p}$\ )$\mathfrak{.\ }$Nevertheless,\ this
automorphism endows $\ V_{n}\ $with the semidirect Lie algebra\ $\mathfrak{%
p\rtimes k}$\ which will be of some relevance \ for the material below.

\subsubsection{\textbf{Coadjoint orbits through rank one \ matrices}.}

Consider  first the symmetric pair $(SL_{n+1}(R),\ SO_{n+1}(R))$\ with 
the quadratic form $\ \langle A,B\rangle =-\frac{1}{2}Tr(AB)$\ on $sl_{n+1}(R)$. This form is positive definite on the space of symmetric
matrices $\fp$ and negative definite on $\fk=so_{n+1}(R)$.  Suppose that $%
P=x\otimes x\ $is a rank one symmetric matrix generated by a\ vector$\ x\in 
\mathbb{R}^{n+1}.$ Then$\ P_{0}=x\otimes x\ -\frac{||x||^{2}}{(n+1)}I\ \ $is\
 in $\fp$ since $Tr(x\otimes x)=\sum_{i=1}^{n+1}(e_{i},\ (x\otimes
x)e_{i})=||x||^{2}.$\ $\ $There are two coadjoint orbits through $P_{0},\ $%
one relative to the action of $Sl_{n+1}(R)\ \ $and the other relative to the action of the
semidirect product\ \ $\fp\rtimes SO_{n+1}(R).\ $\ 

\begin{proposition}
\label{rank one ssimple} The coadjoint orbit $S$  through $%
P_{0}=x_{0}\otimes x_{0}\ -\frac{||x_{0}||^{2}}{(n+1)}I$\ \ is symplectomorphic
to the cotangent bundle of the projective space $\mathbb{P}^{n+1}$ in the semisimple case, and is symplectomorphic to the cotangent bundle of the sphere $S^n$ in the semidirect case.
\end{proposition}

\begin{proof}
Let\ $\mathcal{S\ }$denote\ the coadjoint orbit\ of $Sl_{n+1}(R)\ $through $%
P_{0}=x_{0}\otimes x_{0}\ -\frac{||x_{0}||^{2}}{(n+1)}I.\ $ If $R\in
SO_{n+1}(R)\ $then\ \ $RP_{0}R^{-1}=Rx_{0}\otimes Rx_{0}-\frac{||x_{0}||^{2}%
}{(n+1)}I.$\ It follows that\ $x_{0}$\ can be replaced by $||x_{0}||e_{0},\ $%
because $SO_{n+1}(R)\ $acts transitively on spheres.\ Therefore,\ \ $%
\mathcal{S\ }$is\ \ dieffeomorphic\ to$\ ||x_{0}||^{2}S(e_{0}\ \otimes
e_{0})S^{-1},\ S\in Sl_{n+1}(R).\ \ $Consider now\ the orbit through$\
e_{0}\ \otimes e_{0}.$\ It follows that $\ S(e_{0}\ \otimes
e_{0})S^{-1}=Se_{0}\otimes (S^{T})^{-1}e_{0}\ $for any $S\in SL_{n+1}(R).\ $%
It is easy to verify that \ for any $x\neq 0\ $and any\ $y\ $such that $%
x\cdot y=1\ $there exists a matrix $S\in SL_{n+1}(R)$\ such that $Se_{0}=x$\
and\ $S^{T}y=e_{0}.$\ Hence, $\ S(e_{0}\ \otimes e_{0})S^{-1}=x\otimes y.$
  
  The set of matrices$\{x\otimes y-\frac{(x\cdot y}{n+1}I:x\cdot y=1\}$ can be identified with the set of lines $\{(\alpha x,\frac{1}{\alpha}y): x\cdot y=1\}$ which is symplectomorphic to the tangent bundle of $\mathbb{P}^{n+1}$. The latter is identified with the cotangent bundle via the ambient Euclidean inner product.
  
Consider now the coadjoint orbit relative to the semidirect case.
It follows from (\ref{sdorbit})\ that $\mathcal{S}\ $consists\ of\ matrices\ 
$P=Ad_{h}(P_{0})=h(x_{0})\otimes h(x_{0})-\frac{||x_{0}||^{2}}{(n+1)}I\ \ $%
and $\ \ Q=[Ad_{h}(P_{0}),X]=[x\otimes x,X]=Xx\wedge x\ $since $Q_{0}=0.\ $%
Therefore, $P$\ is a rank one matrix generated by $x=h(x_{0}),\ $with $h\in
SO_{n+1}(R)\ \ $and $Q$\ is a rank two antisymmetric matrix\ $\ x\wedge y\ $%
\ with\ $y=-Xx.\ $\ Since\ $X$\ \ is an arbitrary symmetric matrix of trace
zero\ $y$ can be any point in $\mathbb{R}^{n+1}.$ The  correspondence
 between\ $(x,y)\in \IR^{n+1}\times\IR^{n+1}\rightarrow x\otimes x-\frac{||x||^2}{n+1}I+x\wedge y $ is one to one provided that $x\cdot
y=0. $\ Moreover, $x$\ can be any point of the sphere $||x||=||x_{0}||\ $\
since $SO_{n+1}(R)\ $acts\ transitively \ on spheres by conjugations.
Therefore \ $\mathcal{S\ }$ is identified with points $(x,y)\ $such that $%
||x||=||x_{0}||\ and\ x\cdot y=0$ \ which is the tangent (cotangent) bundle of the
sphere $S^{n} $\, since the two bundles are identified via the Euclidean inner product in $\IR^{n+1}.$
\end{proof}

It may be instructive to show directly that the canonical symplectic form on the cotangent bundle of the sphere coincides with the symplectic form of the coadjoint orbit.

The tangent bundle of the cotangent bundle of the sphere is given by the \ vectors $(x,y,%
\dot{x},\dot{y})\ $in $\mathbb{R}^{4(n+1)}\ $subject to the constraints ~$\
||x||=||x_{0}||,\ x\cdot \dot{x}=0,\ x\cdot y=0$\ and $\dot{x}\cdot y+x\cdot 
\dot{y}=0\ $ and these vectors are identified with matrices $\dot{x}\otimes x+x\otimes \dot{x}+\dot{x}\wedge y+x\wedge\dot{y}$ on the tangent bundle of the coadjoint orbit ( we have omitted the trace factor since it is irrelevant for the these calculations). The  canonical symplectic form on the cotangent bundle of the sphere is given by \[\omega_{(x,y)}((\dot{x}_1,\dot{y}_1),(\dot{x}_2,\dot{y}_2))=\dot{x}_1\dot{y}_2-\dot{x}_2\dot{y}_1
\] 
It follows from above that $\ l\in \mathcal{O}_{G}(l_{0})$\ is
identified with \ $L=x\otimes x+x\wedge y.\ \ $Then  tangent vectors \ $%
v=ad^{\ast }M(l)\ $\ at $l$ are identified  with matrices \ $V\ $via the formula$\
\langle V,U\rangle =\langle L,[M,U]\rangle $\ for all $U\in \fg_s$\ .

An easy calculation shows that $V=[x\otimes x,B]+[x\wedge y,A]+[x\otimes x,A]$   is the tangent vector at $L$ defined by  $M=A+B$ with $ A^T=-A,B^T=B$.  Then,
 $V=\dot{x}\otimes x+x\otimes \dot{x}+\dot{x}\wedge y+x\wedge\dot{y}$ implies that 
 \begin{equation*}
 \dot{x}=Ax,\dot{y}=Ay-Bx
 \end{equation*}

It follows from (\ref{symp}) that the symplectic form on the coadjoint orbit is given by $\omega_L(V_1,V_2)=\langle L,[M_2,M_1]\rangle$. Then,

$\omega
_{x\otimes x+x\wedge y}(v_{1},v_{2})=\langle L,[B_{2},A_{1}]+[A_{2},B_{1}]+[A_{2},A_{1}]\rangle= 
$

\thinspace $\langle x\otimes x,[B_{2},A_{1}]+[A_{2},B_{1}]\rangle+\langle x\wedge
y,[A_{2},A_{1}]\rangle=$

$A_{2}x\cdot B_{1}x-A_{1}x\cdot B_{2}x-(A_{2}x\cdot A_{1}y-A_{1}x\cdot
A_{2}y)=$

$A_{2}x\cdot (A_{1}y-\dot{y}_{1})-A_{1}x\cdot (A_{2}y-\dot{y}%
_{2})-(A_{2}x\cdot A_{1}y-A_{1}x\cdot A_{2}y)=$

$-A_{2}x\cdot \dot{y}_{1}+A_{1}x\cdot \dot{y}_{2}=(\dot{x}_{1}\cdot \dot{y}%
_{2}-\dot{x}_{2}\cdot \dot{y}_{1}).$

Next consider analogous orbits defined by the pseudo Riemannian symmetric
pair $(SL_{n+1}(R),\ SO(1,n))\ $with the Cartan space $\mathfrak{p}$\
consisting of matrices $P$\ $=\left( 
\begin{array}{cc}
0 & -p^T \\ 
p & P_{0}%
\end{array}%
\right) \ $with $p$\ an $n\times 1\ $matrix  and  $P_{0}\ $an $n\times n\ $%
symmetric matrix. These orbits will be defined through the hyperbolic inner product $\
(x,y)_{-1}=x\cdot Jy=x_{1}y_{1}-\sum_{i=2}^{n+1}x_{i}y_{i}\ $where \ $%
J=diag(1,-1,\ldots -1)$.$\ \ $It is easy to verify that $P\in \mathfrak{p}\ $%
if and only if $Tr(P)=0\ $\ and\ $\ (Px,y)_{-1}=(x,Py)_{-1}$; 
similarly, $A\in so(1,n)\ \ $if and only if $(Ax,y)_{-1}=-(x,Ay)_{-1}\ $\ for
all $x,y\ $\ in $\mathbb{R}^{n+1}.\ $Thus $\mathfrak{p}$ is the space of
"hyperbolic symmetric" matrices.

\begin{definition}
$\ \ $The hyperbolic rank one matrices are matrices of the form$\ x\otimes
Jx $\ for some vector \ $x\in \mathbb{R}^{n+1}.\ $\ They will be denoted by $%
(x\otimes x)_{-1}.$
\end{definition}

$\ $It follows that $\ (x\otimes x)_{-1}u=(x,u)_{-1}x$,\ and therefore, 
\begin{equation*} (
 (x\otimes x)_{-1}u,v)_{-1}=(x,u)_{-1}(x,v)_{-1}=(u,\ (x\otimes
x)_{-1}v)_{-1}. \end{equation*} Since the trace of $\ (x\otimes x)_{-1}$\ is equal
to $||x||_{-1}^{2}=(x,x)_{-1}$, $P=\ (x\otimes x)_{-1}$\ $-\frac{%
||x_{0}||_{-1}^{2}}{(n+1)}I\ \ $is in $\mathfrak{p}.$

$\ $Now define rank two skew symmetric\ hyperbolic\ matrices
$Q=(x\wedge y)_{-1}=x\otimes Jy-y\otimes Jx$.\ An easy calculation shows that 
$Q\in so(1,n).$\ Then

\begin{proposition}
\label{hypcot} The coadjiont orbit $\mathcal{S}\ $of the semidirect\ product 
$G_{s}=\emph{Sym}_{-1}\ltimes SO(1,n)\ \ $through $\ P_{0}=(x_{0}\otimes
x_{0})_{-1}\ -\frac{||x_{0}||_{-1}^{2}}{(n+1)|}I\ \ $is equal to $\{(x\otimes x)_{-1}-\frac{||x||_{-1}}{n+1}I+(x\wedge y)_{-1}:||x||_{-1}=||x_0||_{-1}\}$. The latter is symplectomorphic to  the
cotangent bundle of the hyperboloid \ $\mathbb{H}^{n}=%
\{(x,y):||x||_{-1}=||x_{0}||_{-1},(x,y)_{-1}=0\}.$\ The canonical symplectic
form $\omega $\ on $S\ $\ is given by 
\begin{equation*}
\ \omega _{x,y}((\dot{x}_{1},\dot{y}_{1}),(\dot{x}_{2},\dot{y}_{2}))=(\dot{x}%
_{1},\dot{y}_{2})_{-1}-(\dot{x}_{2}\cdot \dot{y}_{1})_{-1}
\end{equation*}%
where $(\dot{x}_{1},\dot{y}_{1})\ $and $(\dot{x}_{2},\dot{y}_{2})\ $denote
tangent vectors at $(x,y).$
\end{proposition}

\begin{proof}
The proof is analogous to the proof in the previous proposition and will be
omitted.
\end{proof}

The coadjoint orbits of the above semidirect products through matrices of
rank one can be expressed in terms of a single parameter $\epsilon \ =\pm 1\ 
$with $SO_{\epsilon }$ $=$ $SO_{n+1}(R)$\ for $\epsilon =1\ $and\ $%
SO_{\epsilon }=SO(1,n)$ for $\epsilon =-1.$\ Then  $\fp_\ep$  will denote the vector space  $\{\left (\begin{matrix}0&\ep p^T\\p&P_0 \end{matrix}\right ):p\in\IR^n, P_0^T=P_0, Tr(P_0)=0\}$. A matrix $X$ belongs to $\fp_\ep$ if and only if it is symmetric relative to the quadratic form $(x,y)_\ep=x_1y_1+\ep\sum_{i=2}^{n+1}x_iy_i$.

We will let $\mathcal{S}_{\epsilon }$\
 denote the coadjoint orbit of the semidirect product 
$G_{\epsilon }=\fp_{\epsilon }\rtimes \ SO_{\epsilon }$ through $%
P_{0}=(x_{0}\otimes x_{0})_{\epsilon }\ -\frac{||x_{0}||_{\epsilon }^{2}}{%
(n+1)}I\ $\ where $||x_{0}||_{\epsilon }^{2}=\ (x\ ,x\ )_{\epsilon }$. The symplectic form $\omega _{x,y}((\dot{x}_{1},%
\dot{y}_{1}),(\dot{x}_{2},\dot{y}_{2}))=(\dot{x}_{2},\dot{y}_{1})_{\epsilon
}-(\dot{x}_{1}\cdot \dot{y}_{2})_{\epsilon }\ $is dual to the Poisson
form 
\begin{equation}
\{f_{1},f_{2}\}_{\epsilon }=(\frac{\partial f_{1}}{\partial x},\frac{%
\partial f_{2}}{\partial y})_{\epsilon }-(\frac{\partial f_{2}}{\partial x},%
\frac{\partial f_{1}}{\partial y})_{\epsilon }  \label{Poiss}
\end{equation}

\section{The\ affine\ Hamiltonian\ system$\ $on coadjoint orbits of rank one}

Consider now the restriction of\ \ $H=\frac{1}{2}\langle L_{\mathfrak{k}}, L_{%
\mathfrak{k}}\rangle+\langle A,L_{\mathfrak{p}}\rangle \ $to the semidirect orbits \ $\mathcal{S}%
_{\epsilon }$ where
$\langle A,B\rangle =-\frac{1}{2}Tr(AB).\ $ This trace form is negative definite on the space of symmetric matrices and positive definite on $so_{n+1}(R) $, but in the pseudo-Riemannian case it is indefinite on both 
 $\fp_{-1} $ and $\fk_{-1}$. 

In what follows it will be convenient to relax the  condition that  $Tr(A)=0$. It is clear that both $A$ and $A-\frac{Tr(A}{n+1}I$ define the same affine Hamiltonian since $\langle I,L_\fp\rangle=0$. Then the
restrictions\ of $L_{\mathfrak{k}}$\ and $L_{\mathfrak{p}}$\ to \emph{S}$%
_{\epsilon }\ $are given by$\ L_{\mathfrak{k}}=(x\wedge y)_{\epsilon
}$ and \ $L_{\mathfrak{p}%
}=(x\otimes x)_{\epsilon }\ -\frac{||x||_{\epsilon }^{2}}{(n+1)}I\ .$\ $\ $%
An easy\ calculation shows that \ the restriction of \ $H\ $is given by 
\begin{equation}
H=\frac{1}{2}||x||_{\epsilon }||y||_{\epsilon }-\frac{1}{2}(Ax,x)_{\epsilon }.
\label{newman}
\end{equation}

\subsection{M\textbf{echanical system of Newmann. }}

If we replace$\ A$ by $-A\ $then the restriction\ \ of the affine 
Hamiltonian \ \ to $S_\ep$ with $\ep=1$ is given by $H=\frac{1}{2}%
(||x||^{2}||y||^{2}+\frac{1}{2}(x,Ax))\ $which coincides with the
Hamiltonian of \ the mechanical problem on the sphere with a quadratic
potential of C. Neumann (\cite{Mos}\ and \cite{Nmn}).
Then $H=\frac{1}{2}(||x||_{-1}^{2}||y||_{-1}^{2}+\frac{1}{2}%
(x,Ax))_{-1}\ $can be considered as the hyperbolic analogue of the problem of Newmann.

 The 
equations of motion\ (\ref{norext}) reduce to 
\begin{equation*}
\frac{d}{dt}(x(t)\wedge y(t))_{\epsilon }=[-A,(x\otimes y)_{\epsilon
}]=(Ax(t)\wedge x(t))_{\epsilon },
\end{equation*} and
\begin{equation*}
\frac{d}{dt}(x(t)\otimes x(t))_{\epsilon }=[(x(t)\wedge y(t)_{\epsilon
},(x(t)\otimes x(t)_{\epsilon }]=||x||_{\epsilon }^{2}(x(t)\otimes
y(t))\epsilon +(y(t)\otimes x(t))_{\epsilon }).
\end{equation*}

An easy calculation shows that the preceding equations are equivalent to
\begin{equation}
\frac{dx(t)}{dt}=||x||_{\epsilon }^{2}y(t),\ \frac{dy}{dt}(t)=-Ax(t)+\frac{1%
}{||x||_{\epsilon }^{2}}(Ax(t),x(t))_{\epsilon }-||y(t)||_{\epsilon
}^{2})x(t).  \label{Newmann}
\end{equation}

Equations ( \ref{Newmann}\ ) with $\epsilon =1\ $and\ $||x||=1\ $form a
point of departure for J. Moser's   book\ on\ integrable\
Hamiltonian\ systems (\cite{Mos}). We will presently show that all
the groundwork for integrability has already been laid out in this paper  in the section on
spectral representations. But first\ let us show that\ equations (\ \ref%
{Newmann})\ can also be\ derived\ in a self contained way from a "mechanical
point of view" through the Maximum Principle of control.\ 

This\ mechanical problem will be phrased as an optimal control problem\ of
minimizing the Lagrangian $\frac{1}{2}\int_{0}^{T}(||u(t)||^{2}-(Ax,x)_{%
\epsilon })dt\ $over the absolutely continuous curves $x(t)\ $\ on an
interval $[0,T]\ $that satisfy $||x||_{\epsilon }=||x_{0}||_{\epsilon }$,\ $%
\frac{dx}{dt}(t)=u(t),\ $and \ also satisfy fixed boundary conditions $%
x(0)=x_{0}$ and $x(T)=x_{1}.$

The constraint $||x||_{\epsilon }=||x_{0}||_{\epsilon }\ $implies \ that $%
(u(t),x(t))_{\epsilon }=0.\ $The Maximum\ Principle of optimal control leads
to the appropriate Hamiltonian on the cotangent bundle of the sphere $%
||x||_{\epsilon }=||x_{0}||_{\epsilon }=1$ (the hyperbolic sphere is the unit hyperboloid $\IH^n$). 

We will use $T^*S_\ep$ to denote this cotangent bundle. It  will be identified with the subset of $\mathbb{R}%
^{n+1}\times \mathbb{R}^{n+1}\ $subject to the constraints 
 $G_{1}=||x||_{\epsilon }-1=0 $ 
and $G_{2}=(y,x)_{\epsilon }=0.$ The Maximum Principle states that the
appropriate Hamiltonian for this problem\ is obtained by maximizing

$H_{0}=-\frac{1}{2}(||u(t)||_{\epsilon }^{2}-(Ax,x)_{\epsilon
})+(y,u)_{\epsilon }\ $relative to the controls $u\ $ that are subject to the constraint $G_0=(x,u)_\ep=0$. According to the method of Lagrange the
maximal Hamiltonian is  obtained by maximizing \emph{H=}$%
H_{0}+\lambda _{0}G_{0}$.
It follows that the optimal control $u\ $\ occurs at $u=y+\lambda_0x$. But then $(u,x)_\ep=(y,x)_\ep+\lambda_0||x||_\ep$ implies that $\lambda_0=0$. Hence, the 
maximal value of $H_{0}\ $is given by
 \begin{equation}\label{Ham}
 H_{0}=\frac{1}{2}(||y||_{\epsilon
}^{2}+(Ax,x)_{\epsilon }).
\end{equation}

The above Hamiltonian is to be taken as a Hamiltonian on $T^*S_\ep$; hence, $T^*S_\ep$ must an invariant manifold for $\vec{H_0}$.  Therefore,  integral curves of $\vec{H_0}$ are the restrictions to $T^*S_\ep$ of the integral curves of a modified Hamiltonian
 \[H=H_0+\lambda_1 G_1+\lambda_2 G_2\] 
 in which the multipliers$\ \lambda _{1}$and $\lambda
_{2} $ are determined by requiring  that   $T^*S_\ep$ be invariant for $\vec{H}$. This requirement will be satisfied whenever the Poisson brackets $\{\emph{H,}%
G_{1}\}\ $and $\{\emph{H,}G_{2}\}\ $  vanish on $T^*S_\ep$. The vanishing of these Poisson brackets implies  that   $\lambda
_{1}=-\frac{\{H_{0,}G_{2}\}}{\{G_{2},G_{1}\}}$ and $\lambda _{2}=\frac{%
\{H_{0,}G_{1}\}}{\{G_{1},G_{2}\}}$.\ \ 

An easy calculation\ based on \ (\ref{Poiss})\ yields \ $\{H_{0},G_{1}%
\}=-2(y,x)_{\epsilon },\ \{H_{0},G_{2}\}=(Ax,x)_{\epsilon }-(y,y)_{\epsilon
}\ $and \ $\{G_{1},G_{2}\}=2||x||_{\epsilon }^{2}=2$ from which it follows that 
\[\lambda _{1}=(Ax,x)_{\epsilon }-(y,y)_{\epsilon } ,\lambda _{2}=(y,x)_{\epsilon }.\]

Hence,%
\begin{equation*}
\emph{H}=\frac{1}{2}||y||_{\epsilon }^{2}+\frac{1}{2}(Ax,x)_{\epsilon }-%
((Ax,x)_{\epsilon }-(y,y)_{\epsilon })G_{1}-%
(y,x)_{\epsilon }G_{2}.
\end{equation*}

The flow\ of $H\ $restricted to $G_{1}=0,G_{2}=0\ $is given by 
\begin{equation}
\frac{dx}{dt}=\frac{\partial \emph{H}}{\partial y}=y,\ \frac{dy}{dt}=-\frac{%
\partial \emph{H}}{\partial x}=-Ax+\frac{1}{||x|_{\epsilon }|^{2}}%
((Ax,x)_{\epsilon }-(y,y)_{\epsilon })x(t)  \label{newmann1}
\end{equation}
The preceding equations coincide with (\ref%
{Newmann}).

\subsection{\textbf{ Integrability}}

The spectral invariants of $L_{\lambda }=L_{\mathfrak{p}}-\lambda L_{%
\mathfrak{k}}-\lambda ^{2}A\ \ $naturally lead to the appropriate
coordinates in terms of which the above equations can be integrated.\ \ It
will be more convenient to divide$\ L_{\lambda }$\ by $-\lambda ^{2}\ $and
redefine $\ L_{\lambda }=-\frac{1}{\lambda ^{2}}L_{\mathfrak{p}}+\frac{1}{%
\lambda }L_{\mathfrak{k}}+A.\ $Since the trace of $(x\otimes x)_\ep $ is a
scalar multiple of the identity it is inessential for the calculations below
and will be omitted. Then\ the\ spectrum\ of\ \ $L_{\lambda \ }$ is given \
by \ \ 
\begin{equation}
0=Det(zI-L_{\lambda })=Det(zI-A)Det(I-(zI-A)^{-1}(-\frac{1}{\lambda ^{2}}L_{%
\mathfrak{p}}+\frac{1}{\lambda }L_{\mathfrak{k}})),
\end{equation}%
\ 

with\ $L_{\mathfrak{p}}=(x\otimes x)_{\epsilon }\ $and\ $L_{\mathfrak{k}%
}=(x\wedge y)_{\epsilon }.$ It\ follows\ that spectral calculations can be
reduced to

\ 
\begin{equation}
0=Det(I-(zI-A)^{-1}(-\frac{1}{\lambda ^{2}}L_{\mathfrak{p}}+\frac{1}{\lambda 
}L_{\mathfrak{k}})),\ Det(zI-A)\neq 0.  \label{spectraleq}
\end{equation}

$~\ $

\ Matrix $I-(zI-A)^{-1}(-\frac{1}{\lambda ^{2}}L_{\mathfrak{p}}+\frac{1}{%
\lambda }L_{\mathfrak{k}})$\ \ is of\ the form$\ $%
\begin{equation}
I-R_{z}(x_{1}\otimes \xi _{1}+x_{2}\otimes \xi _{2})
\end{equation}%
$\ $where\ \ $x_{1}=x,\ x_{2}=y,\ $\ $\ \xi _{1}=-\frac{1}{\lambda ^{2}}%
x_{1}+\frac{1}{\lambda }x_{2},\ \ \xi _{2}=-\frac{1}{\lambda }x_{1\ }$and\
$\ R_{z}=(zI-A)^{-1}.$

A simple argument involving a change of basis shows that the solution of
equation (\ref{spectraleq}\ can be reduced to $\ Det(I-W_{z})=0,$\ 
$\ $where\ $W_{z}\ =(w_{ij})$\ is a $2\times 2\ $ matrix \ with entries $%
w_{ij}$\ \ equal\ \ to the$\ (Rx_{i},\xi _{j})_{\epsilon }\ $(\cite{MosCh}%
).\

 It follows that
\begin{equation}
W_{z}=\left( 
\begin{array}{cc}
(R_{z}x,x)_{\epsilon } & (R_{z}x,y)_{\epsilon } \\ 
(R_{z}y,x)_{\epsilon } & (R_{z}y,y)_{\epsilon }%
\end{array}%
\right) \left( 
\begin{array}{cc}
-\frac{1}{\lambda ^{2}} & \frac{1}{\lambda } \\ 
-\frac{1}{\lambda } & 0%
\end{array}%
\right) 
\end{equation}%
or, 
\begin{equation}
W_{z}=\left( 
\begin{array}{cc}
-\frac{1}{\lambda ^{2}}(R_{z}x,x)_{\epsilon }-\frac{1}{\lambda }%
(R_{z}x,y)_{\epsilon } & \frac{1}{\lambda }(R_{z}x,x)_{\epsilon } \\ 
-\frac{1}{\lambda ^{2}}(R_{z}y,x)_{\epsilon }-\frac{1}{\lambda }%
(R_{z}y,y)_{\epsilon } & \frac{1}{\lambda }(R_{z}y,x)_{\epsilon }%
\end{array}%
\right) .
\end{equation}

Then$\ \ Det(I-W_{z})=0\ $if and only if \ $1-Tr(W_{z})+Det(W_{z})=0$ which
in turn\ \ implies that$\ 1+\frac{1}{\lambda ^{2}}(R_{z}x,x)_{\epsilon }+%
\frac{1}{\lambda ^{2}}((R_{z}x,x)_{\epsilon }(R_{z}y,y)_{\epsilon
}-(R_{z}x,y)_{\epsilon }^{2})=0.$

Let%
\begin{equation}
F=(R_{z}x,x)_{\epsilon }+(R_{z}x,x)_{\epsilon }(R_{z}y,y)_{\epsilon
}-(R_{z}x,y)_{\epsilon }^{2}.  \label{ratint}
\end{equation}

$\ \ $It follows from above \ that $0=Det(zI-L_{\lambda })\ $\ outside of
the spectrum of $A\ $if and only if \ $F(z)=-\lambda ^{2}$.  It is easy to verify that $\lim_{z\rightarrow \pm \infty }zF(z)=(x,x)_{\epsilon
}=1 $  which implies that  $F(z)$ takes both positive and negative values for any $x\neq 0$.
Therefore, $F$ is constant along any solution of  (\ref{Newmann}).

Function $F$\ is rational with poles at the eigenvalues of the matrix $A.$
Hence, $F(z)\ $ will be constant  along the solutions of \ (\ref{Newmann})\ if and only
if the residues of $F$\ are constant along the solutions of\ (\ref{Newmann}).

In\ the  Euclidean\ case the eigenvalues of $A $ are real and distinct
since $A $ is symmetric and regular. Hence, there is no loss in generality
in assuming that \ $A$\ is diagonal.\ Let $\alpha _{1},\ldots ,\alpha _{n+1}$%
\ denote its diagonal entries \ Then,
\begin{equation*}
F(z)=\sum_{k=0}^{n}\frac{F_{k}}{z-\alpha _{k}},
\end{equation*}

where $F_{0},\dots ,F_{n}\ $denote the residues of $F.$ \ It \ follows that  $F_k=\lim_{z\rightarrow\alpha_k}F(z)$.

Since $F(z)=\sum_{k=0}^{n}\frac{x_{k}^{2}}{z-\alpha _{k}}+\sum_{k=0}^{n}%
\sum_{j=0}^{n}\frac{x_{k}^{2}y_{j}^{2}}{(z-\alpha _{k})(z-\alpha _{j})}%
-(\sum_{k=0}^{n}\frac{x_{k}y_{k}}{z-\alpha _{k}})^{2}=$

$\sum_{k=0}^{n}\frac{x_{k}^{2}}{z-\alpha _{k}}+\sum_{k=0}^{n}\sum_{j=0,j\neq
k}^{n}\frac{x_{k}^{2}y_{j}^{2}}{(z-\alpha _{k})(z-\alpha _{j})}%
-2\sum_{k=0}^{n}\sum_{j=0,j\neq k}^{n}\frac{x_{k}y_{k}x_{j}y_{j}}{(z-\alpha
_{k})(z-\alpha _{j})},$

$\lim_{z\rightarrow \alpha _{k}}(z-\alpha
_{k})F(z)=x_{k}^{2}+\sum_{j=0,j\neq k}^{n}\frac{%
x_{j}^{2}y_{k}+x_{k}^{2}y_{j}^{2}}{(\alpha _{k}-\alpha _{j})}%
-2\sum_{j=0,j\neq k}^{n}\frac{x_{k}y_{k}x_{j}y_{j}}{(\alpha _{k}-\alpha _{j})%
}=$

$x_{k}^{2}+\sum_{j=0,j\neq k}^{n}\frac{(x_{j}y_{k}-x_{k}y_{j})^{2}}{(\alpha
_{k}-\alpha _{j})}$. 

Therefore,
\begin{proposition}
\label{Intmotion}Each residue\ \ $F_{k}=x_{k}^{2}+\sum_{j=0,j\neq k}^{n}%
\frac{(x_{j}y_{k}-x_{k}y_{j})^{2}}{(\alpha _{k}-\alpha _{j})},k=0,\ldots ,n\
\ $is an integral of motion for the Hamitonian system \ (\ \ref{Newmann}).\
Moreover, functions $F_{0},\ldots ,F_{n}$\ are in involution relative to the
canonical Poisson bracket in $\mathbb{R}^{n+1}\times \mathbb{R}^{n+1}.$

\begin{proof}
The Poisson bracket relative to the orbit structure coincides with the
canonical Poisson bracket on\ $\mathbb{R}^{n+1}\times \mathbb{R}^{n+1}.$
\end{proof}
\end{proposition}
\begin{remark}  Functions $F_k$ are not functionally independent since $\sum_{k=0}^nF_k=||x||^2=1$.
\end{remark}
In the hyperbolic case the situation is slightly different because $A\ $%
cannot be diagonalized over the\ reals.\ In fact every regular matrix in 
\emph{Sym}$_{-1}$\ is conjugate to $A=\left( 
\begin{array}{ccc}
0 & -\alpha & 0^{T} \\ 
\alpha & 0 & 0 \\ 
0 & 0 & D%
\end{array}%
\right) $, \ $\ $where$\ \alpha $\ is a nonzero number\ and\ $D$\ is a
diagonal $(n-1)\times (n-1)\ $matrix\ with\ distinct nonzero diagonal
entries $\alpha _{2},\ldots ,\alpha _{n}.$

The most convenient way to pass from the Euclidean \ to the hyperbolic case\
is to introduce complex coordinates 
\begin{equation}
v_{0}=\frac{1}{\sqrt{2}}(x_{0}+ix_{1}),v_{1}=\frac{1}{\sqrt{2}}%
(x_{0}-ix_{1}),v_{2}=ix_{2},\ldots ,v_{n}=ix_{n},
\end{equation}

\begin{equation}
w_{0}=\frac{1}{\sqrt{2}}(y_{0}+iy_{1}),w_{1}=\frac{1}{\sqrt{2}}%
(y_{0}-iy_{1}),w_{2}=iy_{2},\ldots ,w_{n}=iy_{n}.
\end{equation}

Then\ 

$((z-A)^{-1}x,y)_{-1}=\frac{1}{z^{2}+\alpha ^{2}}(z(x_{0}y_{0-}x_{1}y_{1})-%
\alpha (x_{0}y_{1+}x_{1}y_{0}))-\sum_{j=2}^{n}\frac{1}{z-\alpha _{j}}%
x_{j}y_{j}=$

$\frac{1}{z-i\alpha }v_{0}w_{0}+\frac{1}{z+i\alpha }v_{1}w_{1}+\sum_{j=2}^{n}%
\frac{1}{z-\alpha _{j}}v_{j}w_{j}=\sum_{j=0}^{n}\frac{1}{z-\alpha _{j}}%
v_{j}w_{j},\ $provided that\ 

$\alpha _{0}=i\alpha $\ and $\alpha _{1}=-i\alpha .\ $

Therefore, $\ $the spectral function  $F(z)$ \ (\ref{ratint})\ is formally   
the same in both the hyperbolic  and the Euclidean case.  It follows that
\begin{equation}
F(z)=\sum_{k=0}^{n}\frac{v_{k}^{2}}{z-\alpha _{k}}+\sum_{k=0}^{n}%
\sum_{j=0}^{n}\frac{v_{k}^{2}w_{j}^{2}}{(z-\alpha _{k})(z-\alpha _{j})}%
-(\sum_{k=0}^{n}\frac{v_{k}w_{k}}{z-\alpha _{k}})^{2}.
\end{equation}

The residues $F_0,\dots,F_n$ defined by  $F(z)=\sum_{k=0}^{n}\frac{F_{k}}{z-\alpha _{k}}$ are given by 
\begin{equation}
F_{k}=v_{k}^{2}+\sum_{j=0,j\neq k}^{n}\frac{(v_{j}w_{k}-v_{k}w_{j})^{2}}{%
(\alpha _{k}-\alpha _{j})},\ k=0,\ldots ,n\ .  \label{hypinvar}
\end{equation}

Since $F(z)$ is real valued for real $z$, $F_1=\bar{F}_0$ and each $F_k,k=2,\dots,n,$ is real.  It follows that   

$Re(F_0), Im(F_0),F_k, k=2,\dots,n$ are integrals of motion for the hyperbolic Newmann problem.

 We leave it to the reader to show that
\begin{equation}
F_{0}=\phi _{0}+\frac{1}{2}\sum_{k=2}^{n}\frac{1}{\alpha ^{2}+\alpha _{k}^{2}%
}((-x_{0}y_{k}-x_{k}y_{0})+i(x_{1}y_{k}-y_{1}x_{k}))^{2}(\alpha _{k}+i\alpha
),
\end{equation}

where$\ \phi _{0}=\frac{1}{2}(x_{0}^{2}-x_{1}^{2})+i(x_{0}x_{1}+\frac{1}{%
2\alpha }(x_{0}y_{1}-x_{1}y_{0})^{2}),\ $ and

\begin{equation}
F_{k}=\phi _{k}-\sum_{j=2,j\neq k}^{n}\frac{(x_{j}y_{k}-x_{k}y_{j})^{2}}{%
(\alpha _{k}-\alpha _{j})},
\end{equation}%
$\ $with $\phi _{k}$\ equal to

$-x_{k}^{2}-\frac{1}{\alpha ^{2}+\alpha _{k}^{2}}\alpha
_{k}((x_{k}y_{0}-x_{0}y_{k})^{2}-(x_{k}y_{1}-x_{1}y_{k})^{2})+2\alpha
((x_{k}y_{0}-x_{0}y_{k})(x_{k}y_{1}-x_{1}y_{k})$ for $k\geq 2$.

\subsection{Integration procedure}

In the Euclidean case\ the integration procedure goes back to C.L. Jacobi \
in connection to the geodesic problem on an ellipsoid. Its modern version is
presented in Moser's papers (\cite{Mos},\ \cite{MosCh}). \ Rather than just
to refer to the classical literature for details, it seems worthwhile to
proceed with the main ingredients of this procedure. For simplicity of
exposition we will confine our attention to the Euclidean sphere; the passage
 to the hyperbolic case requires only minor modifications.

\ The integration is done on the manifold S defined by 
\begin{equation*}
||x||=1,(x,y)=0, F_{k}=x_{k}^{2}+\sum_{j=1,j\neq k}^{n+1}\frac{(x_{j}y_{k}-x_{k}y_{j})^{2}}{%
(\alpha _{k}-\alpha _{j})}=c_{k}, k=1,\dots ,n+1
\end{equation*}
defined by the numbers $c_1,\dots, c_{n+1}$ that satisfy $\sum_{k=1}^{n+1}  c_k=1$.

The following auxiliary lemma will be useful for some calculations below

\begin{lemma}
\label{res}Let \ $\ g(z)=\Pi _{k=1}^{n}(z-x_{k})\ $where$\ x_{1},\ldots
,x_{n}$ are any distinct $n$ numbers. Then

\begin{equation*}
\frac{f(x)}{g(x)}=\sum_{k=1}^{n}\frac{f(x_{k})}{g^{^{\prime
}}(x_{k})(x-x_{k})},\ and\ \lim_{x\rightarrow \infty }x\frac{f(x)}{g(x)}%
=\sum_{k=1}^{n}\frac{f(x_{k})}{g^{^{\prime }}(x_{k})}\ \ 
\end{equation*}%
\ $\ $for any polynomial function $f(z).\ $

Furthermore,$\ \sum_{k=1}^{n}\frac{f(x_{k})}{g^{^{\prime }}(x_{k})}=0$ if $%
\deg (f)<n-1$\ and $\sum_{k=1}^{n}\frac{f(x_{k})}{g^{^{\prime }}(x_{k})}=1\ $%
if $\deg (f)$\ $=n-1$ and its leading coefficient is equal to one.

\begin{proof}
Follows easily from the partial fraction expansion of \ $\frac{f(x)}{g(x)}.$
\end{proof}
\end{lemma}

Following Jacobi,  the integration will be carried out in terms of elliptic coordinates $u_{1},\dots ,u_{n} $ defined as the  zeros
of $(R_{z}x,x)=\sum_{k=1}^{n+1}\frac{x_{k}^{2}}{z-\alpha _{k}}$ for each point
on the sphere $||x||=1$. 
In addition, \ use will be made of the zeros of the rational function $\
\sum_{k=1}^{n+1}\frac{c_{k}}{z-\alpha _{k}}\ .$ These zeros will be denoted by
\ $v_{1},\ldots ,v_{n}\ \ $and will be assumed all distinct.\ \ Let%
\begin{equation}
m(z)=\Pi _{k=1}^{n}(z-u_{k}),\ a(z)=\Pi _{k=1}^{n+1}(z-\alpha _{k}),\ b(z)=\Pi
_{k=1}^{n}(z-v_{k}).
\end{equation}

\begin{lemma}\label{exp}

\begin{equation*}
\sum_{k=1}^{n+1}\frac{x_{k}^{2}}{z-\alpha _{k}}=\frac{m(z)}{a(z)},\ and\
\sum_{k=1}^{n+1}\frac{c_{k}}{z-\alpha _{k}}=\frac{b(z)}{a(z)}.
\end{equation*}
\begin{proof}

The fact that any rational function is determined up to a constant factor by
its zeros and poles implies that \ $\sum_{k=1}^{n+1}\frac{x_{k}^{2}}{z-\alpha
_{k}}=c\frac{m(z)}{a(z)},$ where $c\ $is a\ constant. \ It follows from Lemma %
\ref{res}\ that $c=1,$\ because
\[
1=\sum_{k=0}^{n}x_{k}^{2}=c\sum_{k=0}^{n}\frac{m(\alpha _{k})}{a^{,}(\alpha
_{k})}=c.\ \]  The same argument carries over to $\sum_{k=1}^{n+1}\frac{c_k}{z-\alpha_k}$.
\end{proof}
\end{lemma}
Since $x_{k}^{2}=\frac{m(\alpha _{k})}{a^{\prime }(\alpha _{k})},\ \ x_{k}\ $
can be recovered\ up to a sign from $u_{1},\ldots ,u_{n}.\ \ $

Recall now the function $F(z)=(1+(R_z y,y))(R_z x,x)-(R_z x,y))^2$. Then\ \ $%
F(u_{k})=-(R_{u_{k}}x,y)^{2}=\sum_{k=0}^{n}\frac{c_{k}}{z-\alpha _{k}}=\frac{%
b(u_{k})}{a(u_{k})}$ and therefore, 
\[(R_{u_{k}}x,y)=\pm\sqrt{-\frac{b(u_{k})}{%
a(u_{k})},\ }k=1,\ldots ,n.
\] Each choice of the sign  defines a set of $n\ $linear equations
for the variables $y_{1},\ldots ,y_{n+1},$ which together with $(x,y)=0\ $%
determine $y\ $ uniquely in terms of $u=(u_{1},\ldots ,u_{n})$,  and  each choice identifies 
$(u_{1},\dots ,u_{n})$  as a system of coordinates for the Lagrangian manifold $
S$.  

\begin{proposition}
\label{integrat} \ Vectors $\frac{\partial x}{\partial u_{1}},\frac{\partial
x}{\partial u_{2}},\ldots ,\frac{\partial x}{\partial u_{n+1}}$ \ \ form an
orthogonal frame \ on the sphere$\ ||x||=1\ $,  with
$||\frac{\partial x}{\partial u_{k}}||^{2}=-\frac{1}{4}\frac{m^{\prime }(u_{k})%
}{a(u_{k})}$.\  Suppose that $(R_{u_{k}}x,y)=-\sqrt{-\frac{b(u_{k})}{%
a(u_{k})},\ }k=1,\ldots ,n$. Then, the differential equation $\frac{dx}{dt}=y\ $\  is given on $S$ by 
\begin{equation*}
\frac{1}{2}m^{\prime }(u_{k})\frac{du_{k}}{dt}=\sqrt{-a(u_{k})b(u_{k})}%
,k=1,\ldots ,n.
\end{equation*}
\end{proposition}

The preceding equations can\ also be written as \ \ 
\begin{equation}\label{Jac}
\sum_{k=1}^{n}\frac{u_{k}^{n-j}}{2\sqrt{-a(u_{k})b(u_{k})},}\frac{%
du_{k}}{dt}=\delta _{1j},j=1,\ldots ,n-1.
\end{equation}

\begin{proof}
An easy logarithmic differentiation of $x_{k}^{2}=\frac{m(\alpha _{k})}{%
a^{\prime }(\alpha _{k})}$ yields $\frac{\partial x_{k}}{\partial u_{j}}=-%
\frac{x_{k}}{2(\alpha _{k}-u_{j})}.$ This implies that 
\begin{equation*}
\frac{\partial x}{\partial u_{j}}=\frac{1}{2}(u_{j}I-A)^{-1}x,~j=1,\ldots ,n.
\end{equation*}

Therefore, \ $(\frac{\partial x}{\partial u_{j}},\frac{\partial x}{\partial
u_{k}})=\frac{1}{4}((u_{j}I-A)^{-1}x,(u_{k}I-A)^{-1}x)=$

$\frac{1}{4}((u_{k}I-A)^{-1}(u_{j}I-A)^{-1}x,x)=\frac{-1}{4(u_{j}-u_{k})}%
((u_{j}I-A)^{-1}x-(u_{k}I-A)^{-1}x,x))¼=$

$\frac{-1}{4(u_{j}-u_{k})}((R_{u_{j}}x,x)-(R_{u_{k}}x,x))=0$\ for $k\neq j.\ 
$

For\ $j=k$, 

$\ (\frac{\partial x}{\partial u_{k}},\frac{\partial x}{\partial
u_{k}})=\frac{1}{4}((u_{k}I-A)^{-2}x,x)=-\frac{1}{4}\frac{d}{dz}%
(R_{z}x,x)|_{z=u_{k}}=
-\frac{1}{4}\frac{d}{dz}\frac{m(z)}{a(z)}| _{z=u_{k}}$=

$-\frac{1}{4}%
\frac{m^{\prime }(u_{k})}{a(u_{k})}.\ $\ Since $\frac{m^{\prime }(u_{k})}{%
a(u_{k})}\neq 0$, $\frac{\partial x}{%
\partial u_{1}},\frac{\partial x}{\partial u_{2}},\ldots ,\frac{\partial x}{%
\partial u_{n}}$ form an orthogonal frame on \ $S.$

 Let $P_{1},\dots ,P_{n}\ $%
denote the coordinates of $y$ relative to the frame $\frac{\partial x}{\partial u_{1}},\frac{\partial x}{%
\partial u_{2}},\ldots ,\frac{\partial x}{\partial u_{n}}\ $ for $(x,y)\in S$. It follows that
$(\frac{\partial x}{\partial u_{k}},y)=P_{k}||\frac{\partial x}{\partial
u_{k}}||^{2}.\ $ But\ $(\frac{\partial x}{\partial u_{k}},y)=\frac{1}{2}%
((u_{k}I-A)^{-1}x,y)=\frac{1}{2}(R_{u_{k}}x,y)=-\frac{1}{2}\sqrt{-\frac{%
b(u_{k})}{a(u_{k})}.\ }$

Therefore,$\ $%
\begin{equation*}
P_{k}=\frac{2}{m^{\prime }(u_{k})}\sqrt{-a(u_{k})b(u_{k})}.
\end{equation*}

Suppose now that $(x(t),y(t))$ is a curve in $S\ $with $\frac{dx}{dt}=y.$
 Then,\ $\frac{dx}{dt}=\sum_{k=1}^{n}\frac{\partial x}{\partial u_{k}}\frac{%
du_{k}}{dt}=\sum_{k=1}^{n}P_{k}\frac{\partial x}{\partial u_k}.\ $Therefore, \ $\frac{1}{2%
}m^{\prime }(u_{k})\frac{du_{k}}{dt}=\sqrt{-a(u_{k})b(u_{k})}.$

Second expression follows from Lemma \ref{res} which implies that $\
\sum_{k=1}^{n}\frac{u_{k}^{m}}{m^{^{\prime }}(u_{k})}=\delta _{n-1,m},\
m\leq n-1$ by taking $f(z)=z^{m}$\ and $g(z)=m(z).\ $Then, $\ m^{\prime
}(u_{k})=2\sqrt{-a(u_{k})b(u_{k})}/\frac{du_{k}}{dt}$\ which, after the
substitution, leads to

$\sum_{k=1}^{n}\frac{u_{k}^{n-j}}{2\sqrt{-a(u_{k})b(u_{k})},}\frac{%
du_{k}}{dt}\ =\delta _{1j},j=1,\ldots ,n-1.$
\end{proof}
J. Moser points out that equation (\ref{Jac}) is related to the Jacobi map of the Riemann surface 
\begin{equation}
w^2=-4a(z)b(z).
\end{equation}
 In fact he shows that the Jacobi map  given by 
 \[
 \sum_{k=1}^n\int_{(0,0)}^{(u_k,w_k)}\frac{z^{n-j}dz}{2\sqrt{-a(z)b(z)}}
 \]
 takes the divisor class defined by $(u_k, 2\sqrt{a(u_k)b(u_k)}),k=1,\dots,n$ into a point $s\in C^n/ \Gamma$ where $\Gamma$ denotes the period lattice of the differentials of the first kind (\cite{Mos}).
\section{ Connection to geodesic problems  on quadric surfaces:\ Knorrer's
transformation}

Jacobi's geodesic problem on an  ellipsoid   $\mathbb{S}=\{x\in \mathbb{R}^{n+1}:(x,A^{-1}x)=1\}$\  consists of finding curves    in $\IS$ of minimal length,   relative  to the metric inherited from the ambient Euclidean metric in $\IR^{n+1}$, that connect  a given pair of  points in $\IS$.  Jacobi was  able to show that the curves of minimal length  can be obtained from the solutions of a first order partial differential equation, known today as the Hamilton-Jacobi equation; in the process, he discovered an ingenious choice of coordinates on the ellipsoid, known today as elliptic coordinates, in terms of which the associated partial differential equation becomes separable  with its solutions given by hyperelliptic functions. 

Alternatively, the geodesic equations can  be  represented by a Hamiltonian system
\begin{equation}
\frac{dx}{dt}=p,\ \frac{dp}{dt}=-\frac{(p,A^{-1}p)}{||A^{-1}x||^{2}}A^{-1}x.
\label{jac}
\end{equation}
on the cotangent bundle of $\IS$  realized as  the subset of $\IR^{n+1}\times\IR^{n+1}$ subject to 
$G_{1}=(x,A^{-1}x)-1=0,G_{2}=(p,A^{-1}x)=0$.  Moser shows that the above equations are 
generated by the Hamiltonian
\begin{equation*}
\emph{H}=\frac{1}{2}||p||^{2}+\frac{(p,A^{-1}p)}{2||A^{-1}x||^{2}}G_{1}-%
\frac{(p,A^{-1}x)}{||A^{-1}x||^{2}}G_2%
\end{equation*}%
in $\IR^{n+1}\times\IR^{n+1}$
in the sense that,
\[\frac{dx}{dt}=\frac{\partial H}{\partial p},\frac{dp}{dt}=-\frac{\partial H}{\partial x}\]
but constrained to
 $G_1=G_2=0$ and  $H=\frac{1}{2}$, i.e., to  $||p||=1$
 (\cite{Mos}).
 
 Let us modify  above equations by replacing the Euclidean inner product  $(x,y)$ by  the inner product $(x,y)_\ep$ that encompasses both the Euclidean and the hyperbolic inner product. Then, equations (\ref{jac}) take on the following form:
 \begin{equation}
\frac{dx}{dt}=p,\ \frac{dp}{dt}=-\frac{(p,A^{-1}p)_\ep}{||A^{-1}x||_\ep^{2}}A^{-1}x.
\label{jac1}
\end{equation}
assuming that $||A^{-1}x||\neq 0$.
  It follows that  $G_{1}=(x,A^{-1}x)_\ep-1=0,G_{2}=(p,A^{-1}x)_\ep=0,||p||_\ep=1$ is an invariant set  for (\ref{jac1}).

Remarkably, equations (\ref{jac1}) can be transformed into the equations of$\ 
$(\ref{newmann1}) by a transformation discovered by H. Knorrer in (\cite{knor}), and the
integrals of motion \ of the geodesic problem can be deduced from the
integrals of motion associated with the mechanical problem on the sphere. In
what follows we will consider the inverse of the Knorrer's transformation
and show that the integrals of motion for the geodesic problem can be
deduced from \ the mechanical problem of Newmann. For that reason we will
begin with equations (\ref{newmann1}) \ written as

\begin{equation}
\frac{du}{ds}=v,\ \frac{dv}{ds}=-A^{-1}u+((A^{-1}u,u)_{\epsilon
}-||v||_{\epsilon }^{2})u,\ ||u||_{\epsilon }=1,\ \epsilon =\pm 1.
\label{newman2}
\end{equation}
It is important to keep in mind that the matrix $A$ also depends on $\ep$ since it belongs to $\fp_\ep$ (modulo the trace).
In what follows it will be necessary to assume that $(Au,u)_\ep>0 ,u\neq 0$.

Recall that$\ F_{0}=((v,Av)_{\epsilon }-1)(u,Au)_{\epsilon
}-(u,Av)_{\epsilon }^{2}\ $is an integral of motion for (\ref{newman2}) as
can be easily seen from \ (\ref{ratint}) with $z=0.\ $Let $\Phi (\lambda
,u,v)=(x,p)\ $denote the mapping \ from the manifold $N_{0}=\{(\lambda
,u,v):||u||_{\epsilon }=1,(u,v)_{\epsilon }=0,F_{0}=0,\lambda \in \mathbb{R}%
\}$ given by 
\begin{equation}
x=\frac{Au}{\sqrt{(Au,u)_{\epsilon }}},\ p=\frac{\lambda }{\sqrt{%
(Au,u)_{\epsilon }}}(Av-\frac{(Au,v)_{\epsilon }}{(Au,u)_{\epsilon }}Au).
\label{knor}
\end{equation}

\ It follows that $(x,A^{-1}x)_{\epsilon }=1\ $and that \ $%
(p,A^{-1}x)_{\epsilon }=0,\ $hence $\Phi \ $maps into the tangent bundle of
the quadric $(x,A^{-1}x)_{\epsilon }=1\ .$\ We will show now that \ $\Phi \ $%
is invertible and that its inverse is given by: $\ $%
\begin{equation}
\lambda =\frac{\sqrt{(A^{-1}p,p)_{\epsilon }}}{||A^{-1}x||_{\epsilon }},\ u=%
\frac{A^{-1}x}{||A^{-1}x||_{\epsilon }},\ v=\frac{1}{\sqrt{%
(A^{-1}p,p)_{\epsilon }}}(A^{-1}p-(u,A^{-1}p)_{\epsilon }u).  \label{inverse}
\end{equation}

It follows from (\ref{knor}) that $\frac{(Au,u)_{\epsilon }}{\lambda ^{2}}%
(A^{-1}p,p)_{\epsilon }=(Av-\frac{(Au,v)_{\epsilon }}{(Au,u)_{\epsilon }}Au)$
$(v-\frac{(Au,v)_{\epsilon }}{(Au,u)_{\epsilon }}u)=$\ $(v,Av)_{\epsilon }-%
\frac{(Au,v)_{\epsilon }^{2}}{(Au,u)_{\epsilon }}.$

The constraint $((v,Av)_{\epsilon }(u,Au)_{\epsilon }-(u,Av)_{\epsilon
}^{2}=(u,Au)_{\epsilon }\ \ $implies that\ \ $\frac{(Au,u)_{\epsilon }}{%
\lambda ^{2}}(A^{-1}p,p)_{\epsilon }=1.$ \ Further constraints $%
||u||_{\epsilon }^{2}=1$\ and $(u,v)_{\epsilon }=0$\ imply that  $%
(Au,u)_{\epsilon }=\frac{1}{||A^{-1}x||_{\epsilon }^{2}}$\ and \ $\frac{%
(u,A^{-1}p)_{\epsilon }}{\sqrt{(A^{-1}p,p)_{\epsilon }}}=-\frac{%
(Au,v)_{\epsilon }}{(Au,u)_{\epsilon }},\ $hence (\ref{inverse}).

Let $u(s)$ and $v(s)$ be any solutions of (\ref{newman2}) and let $\lambda
(s)\ $be a solution of 
\begin{equation}
\frac{d\lambda }{ds}=2\frac{(Au(s),v(s))_{\epsilon }}{(Au(s),u(s))_{\epsilon
}}\lambda (s).  \label{lambd}
\end{equation}

Then,

$\frac{dx}{ds}=\frac{Av}{\sqrt{(Au,u)_{\epsilon }}}-\frac{(Au,v)_{\epsilon }%
}{\sqrt[3]{(Au,u)_{\epsilon }}}Au=\frac{1}{\lambda }v$, and

$\frac{dp}{ds}=(\frac{1}{\sqrt{(Au,u)_{\epsilon }}}\frac{d\lambda }{ds}%
-\lambda \frac{(Au,v)_{\epsilon }}{\sqrt[3]{(Au,u)_{\epsilon }}})(Av-\frac{%
(Au,v)_{\epsilon }}{(Au,u)_{\epsilon }}Au)+\frac{\lambda }{\sqrt{%
(Au,u)_{\epsilon }}}\frac{d}{ds}(Av-\frac{(Au,v)_{\epsilon }}{%
(Au,u)_{\epsilon }}Au)=\frac{(Au,v)_{\epsilon }}{\sqrt[3]{(Au,u)_{\epsilon }}%
}(Av-\frac{(Au,v)_{\epsilon }}{(Au,u)_{\epsilon }}Au)+\frac{\lambda }{\sqrt{%
(Au,u)_{\epsilon }}}\frac{d}{ds}(Av-\frac{(Au,v)_{\epsilon }}{%
(Au,u)_{\epsilon }}Au)=$

$\frac{(Au,v)_{\epsilon }}{\sqrt[3]{(Au,u)_{\epsilon }}}\lambda (Av-\frac{%
(Au,v)_{\epsilon }}{(Au,u)_{\epsilon }}Au)+\frac{\lambda }{\sqrt{%
(Au,u)_{\epsilon }}}(A\frac{dv}{ds}-\frac{(Au,v)_{\epsilon }}{%
(Au,u)_{\epsilon }}Av-(\frac{1}{(Au,u)_{\epsilon }}((Av,v)_{\epsilon }+(Au,%
\frac{dv}{ds})_{\epsilon })-2\frac{(Au,v)_{\epsilon }^{2}}{(Au,u)_{\epsilon
}^{2}})Au)=$

$-\frac{\lambda }{\sqrt{(Au,u)_{\epsilon }}}u+\frac{\lambda }{\sqrt{%
(Au,u)_{\epsilon }}}Au(\frac{(Au,v)_{\epsilon }^{2}-((Av,v)_{\epsilon
}-1)(Au,u)_{\epsilon }}{(Au,u)_{\epsilon }^{2}})=-\frac{\lambda }{\sqrt{%
(Au,u)_{\epsilon }}}u.$

Hence,%
\begin{equation}
\frac{dx}{ds}=\frac{1}{\lambda }v,\ \frac{dp}{ds}=-\lambda A^{-1}x.
\label{ell1}
\end{equation}

Now identify $\frac{1}{\lambda (s)}\ \ $with a function \ $\frac{dt}{ds}.$
It follows from (\ref{inverse}) that $\frac{dt}{ds}^{2}=\frac{1}{\lambda ^{2}%
}=\frac{||A^{-1}x||_{\epsilon }^{2}}{(A^{-1}p,p)_{\epsilon }}.$ Hence, \ $%
\frac{dt}{ds}\frac{(A^{-1}p,p)_{\epsilon }}{||A^{-1}x||_{\epsilon }^{2}}%
=\lambda (s)$ and hence equations (\ref{ell1}) \ coincide with equations (%
\ref{jac1}) after the reparametrization $t=\phi (s)\ $\ with $\frac{d\phi }{ds%
}=\frac{1}{\lambda (s)}.$

The integrals of motion obtained for the mechanical problem of Newmann have
their analogues for the problem of Jacobi via the following proposition

\begin{proposition}
Let  $F(w)=(1+(R_{w}v,v)_\ep)(R_{w}u,u)_\ep-(R_{w}u,v)_\ep^{2}$, with $R_w=(wI-A)^{-1}$. Then,

\[ 
F(\frac{1}{z})=\frac{1}{||A^{-1}x||_\ep^{2}(A^{-1}p,p)_{\epsilon }}%
(1+(S_{z}x,x))_\ep)(S_{z}p,p)_\ep-(S_{z}x,p)_\ep^{2},
 S_{z}=(z-A)^{-1},\] under the substitutions   given by formulas (\ref{inverse}).

\begin{proof}
\ It is easy to verify that $F(w)$ is invariant under the  change of variable $v\rightarrow v+\alpha u$ with $\alpha$ an arbitrary scalar. Hence, $F(w)=(1+(R_{w}V,V))(R_{w}u,u)-(R_{w}u,V)^{2},\ $%
where $V=\frac{A^{-1}p}{\lambda ||A^{-1}x||}.\ $ In the proof below we will
use the identity
 \[(\frac{1}{z}-A^{-1})^{-1}+A=-(z-A)^{-1}A^{2}.
 \] Then,
\[ 
\begin{array}{lll} &1+(R_{\frac{1}{z}}V,V)_\ep=1+(\frac{1}{z}-A)^{-1})^{-1}V,V)_\ep=\\&
1-((z-A)^{-1}A^{2}\frac{A^{-1}p}{\lambda ||A^{-1}x||_\ep},\frac{%
A^{-1}p}{\lambda ||A^{-1}x||_\ep})_\ep-(\frac{p}{\lambda ||A^{-1}x||_\ep},\frac{A^{-1}p}{%
\lambda ||A^{-1}x||_\ep})_\ep=\\&
-((z-A)^{-1}\frac{Ap}{\lambda ||A^{-1}x||_\ep},\frac{A^{-1}p}{\lambda
||A^{-1}x||_\ep})_\ep=-\frac{1}{\lambda ^{2}||A^{-1}x||_\ep^{2}}(S_{z}p,p)_\ep.
\end{array}
\]
Further,
\[\begin{array}{ll}&(R_{\frac{1}{z}}u,u)_\ep=-((z-A)^{-1}A^{2}\frac{A^{-1}x}{||A^{-1}x||_\ep},\frac{A^{-1}x}{%
||A^{-1}x||_\ep})_\ep-(\frac{x}{||A^{-1}x(s)||_\ep},\frac{A^{-1}x}{||A^{-1}x||_\ep})_\ep\\&
-\frac{1}{||A^{-1}x||_\ep^{2}}(1+(S_z x,x))_\ep,\end{array}
\] and
\[
\begin{array}{ll}&(R_{\frac{1}{z}}u,V)_\ep=
-((z-A)^{-1}A^{2}\frac{A^{-1}x}{||A^{-1}x||_\ep},\frac{A^{-1}p)}{%
\lambda ||A^{-1}x||_\ep})_\ep-(\frac{x}{||A^{-1}x||_\ep},\frac{A^{-1}p}{\lambda
||A^{-1}x||_\ep})_\ep=\\&\frac{1}{\lambda ||A^{-1}x||_\ep^{2}}(S_{z}x,p).
\end{array}
\]
Hence,

$F(\frac{1}{z})=\frac{1}{\lambda ^{2}||A^{-1}x||_\ep^{4}}%
((1+(S_{z}x,x))_\ep(S_{z}p,p)_\ep-(S_{z}x,p)_\ep^{2})=\frac{1}{%
||A^{-1}x||_\ep^{2}((A^{-1}p,p)_{\epsilon }}%
((1+(S_{z}x,x)_\ep)(S_{z}p,p)_\ep-(S_{z}x,p)_\ep^{2}).$
\end{proof}
\end{proposition}

\begin{corollary}
Function$\ G(z)=(1+(S_{z}x,x)_\ep)(S_{z}p,p)_\ep-(S_{z}x,p)_\ep^{2}\ \ is$ constant
along the geodesic flow (\ref{jac1}).
\end{corollary}

\begin{proof}
$||A^{-1}x||_\ep^{2}(A^{-1}p,p)_\ep\ $is an integral of motion for (\ref{jac1})
because
 \[
\begin{array}{ll}&\frac{d}{dt}||A^{-1}x||_\ep^{2}(A^{-1}p,p)_\ep=\\&

2(A^{-1}p,A^{-1}x)_\ep(A^{-1}p,p)_\ep-||A^{-1}x||_\ep^{2}(A^{-1}p,(A^{-1}p,p)_\ep\frac{%
A^{-1}x}{||A^{-1}x||_\ep^{2}})_\ep=0.
\end{array}
\]
 It follows that $G(z)\ $is constant along the
\ solutions of (\ref{jac1}). since $F(\frac{1}{z})\ $is constant along the solutions of
(\ref{Newmann}).
\end{proof}

\begin{remark}
Function $||A^{-1}x||_\ep^{2}(A^{-1}p,p)_{\epsilon }$ is known as Joachimsthal's
integral of motion( \cite{joach}), (\cite{per}).
\end{remark}

In the \ Euclidean case the matrix $A\ $can be assumed diagonal with \ $%
\alpha _{1},\ldots ,\alpha _{n+1}\ $its eigenvalues. An argument identical
to the one used above shows that 
\begin{equation*}
G(z)=\sum_{k=1}^{n+1}\frac{G_{k}}{z-\alpha _{k}}
\end{equation*}

\ and that the residues\ $G_{k}\ $are given by 
\begin{equation}
G_{k}=p_{k}^{2}+\sum_{j=1,j\neq k}^{n+1}\frac{(x_{j}p_{k}-x_{k}p_{j})^{2}}{%
(\alpha _{k}-\alpha _{j})},k=1,\ldots ,n+1,  \label{resjac}
\end{equation}
as reported  in (\cite{Mos}. The hyperbolic case differs only in minor details due to different canonical  structures of $A$.

\section{The case $A=0$ and the problem of Kepler}

Consider now the Hamiltonian \ $H=\frac{1}{2}\langle L_\fk,L_\fk\rangle\ \ $ on coadjoint orbits of  $\fp_\ep\rtimes \fk_\ep$, $\ep=\pm1$, through rank one matrices in $\fp_\ep$. We will continue with the notations of the last two sections and consider the coadjoint orbis through  matrices 
$P_{0}=(x_{0}\otimes x_{0})_{\varepsilon }\ -\frac{||x_0||_\ep^2}{n+1}I,\ep=\pm1\ $.
We have seen that the these coadjoint orbits consist of matrices \[
L_\fp=(x\otimes x)_\ep, L_\fk=(x\wedge y)_\ep,||x||_\ep=||x_0||_\ep, (x,y)_\ep=0,\] that are symplectomorphic to the cotangent bundle of the "sphere"  $\{(x,y):||x||_\ep=||x_0||_\ep,(x,y)_\ep=0\}$.
\ The Hamiltonian equations \ (\ref{norext}) and (\ref{Newmann}) reduce to 
\begin{equation}\label{wedge}
\frac{dL_\ep}{dt}=[L_\fk,L_\fp],\ \frac{dL_\fk}{dt}=0,
\end{equation}
and
\begin{equation}\label{geod}
\frac{dx}{dt}=||x||_{\varepsilon }^{2}y,\ \frac{dy}{dt}=-||y||_{\varepsilon
}^{2}x.
\end{equation}
It follows that the solutions  satisfy
\[\frac{d^{2}x}{dt^{2}}%
+||x||_{\varepsilon }^{2}||y||_{\varepsilon }^{2}x=0
\]
The restriction of the Hamiltonian $H$ to these  orbits is given by
$H=\frac{1 }{2}\ ||x||_{\varepsilon }^{2}||y||_{\varepsilon
}^{2}. $ Hence, on energy level $H=\frac{\ep}{2}$, the preceding equations reduce to
 \begin{equation}\label{geod1}\frac{d^2 x}{dt^2}+\ep x=0.
 \end{equation}
  It follows that  the solutions of (\ref{geod1}) are given by  great circles
$x(t)=a\cos(t)+b\sin t$ for  $\varepsilon =1$,  and  great hyperbolas
$ x(t)=a\cosh t+b\sinh t$ for $ \varepsilon =-1$, with
$||a||_{\varepsilon }^{2}+||b||_{\varepsilon }^{2}=||x_{0}||_{%
\varepsilon }^{2},\ (a,b)\varepsilon =0.$

Recall now \ the Hamiltonian $E=\frac{1}{2}||p||-\frac{1}{||q||}$\ \
associated with the problem of Kepler in the phase space$\ \{(q,p)\in 
\mathbb{R}^{n}\times \mathbb{R}^{n}:\ q\neq 0\}\ $\ corresponding to the
normalized constants $m=kM=1$ and the associated equations of motion
\begin{equation}
\frac{dq}{dt}=p,\ \frac{dp}{dt}=-\frac{1}{||q||^{3}}q.  \label{kepler}
\end{equation}
Below is a summary of  the classical theory connected with Kepler's problem.

\begin{enumerate}
\item $L=q\wedge p\ $and \ $F=Lp-\frac{q}{||q||}\ $\ are constants of motion
for (\ref{kepler}). \ $L$\ is  an $n$ dimensional generalization of the
angular momentum $p\times q$. Its constancy implies that each solution
remains in the plane spanned by $q(0)$ and $p(0).$

The second vector is call the Runge-Lenz vector or, sometimes, the eccentricity
vector. It lies in the plane spanned by $p$\ and $%
q.$

\item
Let $||F||^{2}=2||L||^{2}E+1$, where $||L||^{2}=-\frac{1}{2}Tr(L^{2})=||q||^{2}||p||^{2}-(q\cdot p)^{2}$. Then,
$||F||<1\ $whenever \ $E<0,\,||F||=1\ $%
whenever$\ E=0,$\ and$\ \ ||F||>1\ $whenever $E>0.$

\item A solution\ $(q(t),p(t))$ \ evolves on a line through the origin if
and only if $q(0)$\ and \ $p(0)\ $are colinear, that is, whenever $L=0.\ $%
In the case that $L\neq 0\ $\begin{equation*}
||q(t)||=\frac{||L||}{1+||F||\cos \phi (t)},
\end{equation*}
where $\phi (t)\ $denotes the angle between $F$ and $q(t).$ 
Therefore, 

$q(t)\ $traces an ellipse when $||F||<1$, a parabola when $%
||F||=1$ and a hyperbola when $||F||>1.$
\end{enumerate}

There is a remarkable connection between the solutions of Kepler's problem and the  geodesic flows on space forms that was first  reported by V.A. Fock  in 1935 in connection with the theory of hydrogen atom (\cite{Fk})  which then was rediscovered independently by J. Moser in (\cite{MosKep})  for the geodesics on a sphere. Moser's  study was  later completed to all space forms by Y. Osipov in (\cite{Os}). 

As brilliant as these contributions were, they, nevertheless, did not  attempt any explanations in regard to this enigmatic connection between planetary motions and geodesics on space forms.   This  issue later inspired  V. Guillemin and S. Sternberg  to take up the problem of Kepler in (\cite{GSt})  in a larger geometric context with Moser's  observation  at the heart of the matter.

 It seems altogether natural to include  Kepler's problem in this study. In this setting  Kepler's system is recognized within a large class of integrable systems and secondly, the  focus on coadjoint representations provides  natural explanations for its connections to the geodesic problems. Following Moser  we will 
consider the stereographic projection from the sphere $||x||_{%
\varepsilon }^{2}=h^{2}\ $into $\mathbb{R}^{n}\ $given by 
\[\lambda(x-he_{0})+he_{0}=(0,p)\text{, where }\lambda =\frac{h}{h-x_{0}}.
\]
Here, $(x_0,x_1,\dots,x_n)$ denote the coordinates of a point $x$ in  $\IR^{n+1}$  corresponding to the standard basis $e_0,\dots,e_n$. 
It follows that 
\begin{equation}
x_{0}=\frac{h(||p||^{2}-\varepsilon h^{2})}{||p||^{2}+\varepsilon h^{2}},\ 
\text{and }\bar{x}=x-x_{0}e_{0}=\frac{2\varepsilon h^{2}}{%
||p||^{2}+\varepsilon h^{2}}p.  \label{steregraphic}
\end{equation}

Consider now the extension of this mapping to the  cotangent bundle of 
$||x||_{\varepsilon }^{2}=h^{2}$  that \ pulls back the canonical symplectic
form in $\mathbb{R}^{n}\times \mathbb{R}^{n}\ $onto the symplectic form of
the cotangent bundle of $||x||_{\varepsilon }^{2}=h^{2}.\ $It suffices to
find a mapping $q=\Psi (x,y)$\ such that 
\begin{equation}
\sum_{i=1}^{n}q_{i}dp_{i}=(y,dx)_\ep= y_{0}dx_{0}+\varepsilon \sum_{i=1}^{n}y_{i}dx_{i}
\label{differentls}
\end{equation}

because the symplectic forms are the exterior derivatives of the preceding
forms. \ It turns out that such an extension is unique by the following
arguments.

Let\ \ $x=\Phi (p)\ $denote the mapping \ given by (\ref{steregraphic}\ ).
Then, 
\begin{equation}
dx=(\frac{\partial \Phi }{\partial p})_\ep dp=(\frac{4\varepsilon h^{3}}{%
(||p||^{2}+\varepsilon h^{2})^{2}}p\cdot dp,\frac{2\varepsilon h^{2}}{%
||p||^{2}+\varepsilon h^{2}}dp-\frac{4\varepsilon h^{2}p\cdot dp}{%
(||p||^{2}+\varepsilon h^{2})^{2}}p).  \label{diff1}
\end{equation}

It follows by an easy calculation that
\begin{equation}
||dx||_{\varepsilon }^{2}=dx_0^2+\ep\sum_{i=1}^{n} dx_i^2= \frac{4h^{4}\varepsilon }{(||p||^{2}+\varepsilon
h^{2})^{2}}||dp||^{2}.  \label{normdiff}
\end{equation}

Since
$( dx,dx)_\ep=( (\frac{\partial \Phi }{\partial p})_\ep dp,(\frac{\partial \Phi }{\partial p})_\ep dp )_\ep =( (\frac{\partial \Phi }{\partial p})^*_\ep (\frac{\partial \Phi }{\partial p})_\ep dp,dp)_\ep=\frac{4h^2\ep}{(||p||^2+\ep h^2)}||dp||^2
$
it follows that
 $\ (\frac{\partial \Phi }{\partial p}%
)_{\varepsilon }^{\ast }\frac{\partial \Phi }{\partial p}=\frac{%
4h^{4}\varepsilon }{(||p||^{2}+\varepsilon h^{2})^{2}}I_{n},$ where $I_{n}\ $%
denotes the $n$  dimensional identity  and  $\ (\frac{\partial \Phi }{\partial p}%
)_{\varepsilon }^{\ast }$\ the adjoint operator \ of $\ \frac{\partial \Phi }{%
\partial p}\ $relative to the inner product $\ (\ ,\ )_{\varepsilon },i.e.,(%
\frac{\partial \Phi }{\partial p})_{\varepsilon }^{\ast }=(\frac{\partial
\Phi }{\partial p})^{T}J_{\varepsilon },$ with $(\frac{\partial \Phi }{%
\partial p})^{T}\ $the transpose of $\frac{\partial \Phi }{\partial p}$ and $%
J_{\varepsilon }$ diagonal matrix with its diagonal entries $(1,\varepsilon ,\varepsilon ,\ldots ,\varepsilon ).$

Then, 

$\ q\cdot
dp=(y\cdot dx)_{\varepsilon }=(y ,\frac{\partial \Phi }{\partial p}dp)_\ep=(%
\frac{\partial \Phi }{\partial p})_{\varepsilon }^{\ast }y, dp)_\ep$ implies that $q=(\frac{\partial \Phi }{\partial p})_{\varepsilon }^{\ast }y\ 
$ or,
\[y=\ep\frac{(||p||^{2}+\varepsilon h^{2})^{2}}{4h^{4}}%
 (\frac{\partial \Phi }{\partial p})_\ep q ,\] since ($\frac{\partial \Phi 
}{\partial p})_{\varepsilon }^{\ast }y=\varepsilon \frac{(||p||^{2}+%
\varepsilon h^{2})^{2}}{4h^{4}}(\frac{\partial \Phi }{\partial p}%
)_{\varepsilon }^{\ast }(\frac{\partial \Phi }{\partial p})q=q.$

Equations (\ref{diff1}) \ reveal that 
\begin{equation}
y=(\frac{1}{h}q\cdot p,\frac{||p||^{2}+\varepsilon h^{2}}{2h^{2}}q-\frac{%
q\cdot p}{h^{2}}p)  \label{Hamlift}
\end{equation}

from which \ it follows that 
\begin{equation}
||y||_{\varepsilon }^{2}=\varepsilon \frac{(||p||^{2}+\varepsilon h^{2})^{2}%
}{4h^{4}}||q||^{2}.  \label{geodesiclink}
\end{equation}

To pass to the problem of Kepler, write the Hamiltonian\ $H=\frac{1}{2}%
||x||_{\varepsilon }^{2}||y||_{\varepsilon }^{2}\ \ $in the
variables $(p,q).\ $It follows that $H=\frac{1}{2}h^{2}\varepsilon \frac{%
(||p||^{2}+\varepsilon h^{2})^{2}}{4h^{4}}||q||^{2}=\frac{1}{2}\varepsilon 
\frac{(||p||^{2}+\varepsilon h^{2})^{2}}{4h^{2}}||q||^{2}.$

The corresponding flow is given by 
\begin{equation}
\frac{dp}{ds}=\frac{\partial H}{\partial q}=\varepsilon \frac{%
(||p||^{2}+\varepsilon h^{2})^{2}}{4h^{2}}q,\frac{dq}{ds}=-\frac{\partial H}{%
\partial p}=-\varepsilon \frac{||p||^{2}+\varepsilon h^{2}}{2h^{2}}||q||^{2}p
\label{before par}
\end{equation}
\ 

On energy level $H=\frac{\varepsilon }{2h^{2}},\frac{(||p||^{2}+\varepsilon
h^{2})^{2}}{4}||q||^{2}=1\ $ and the preceding equations reduce to \begin{equation}\frac{dp}{%
ds}=\varepsilon \frac{q}{h^{2}||q||^{2}},\ \frac{dq}{ds}=-\varepsilon \frac{%
||q||}{h^{2}}p.\end{equation} The preceding equations coincide with the equations of Kepler's problem  (\ref{kepler}) after the reparametrization by a parameter  $\ t=-\frac{-\varepsilon }{h^{2}}%
\int_{0}^{s}||q(\tau ||d\tau .$  For then, $\frac{ds}{dt}=-\frac{%
\ep h^2}{||q||}$ and equations (\ref{before par}) become 
\begin{equation*}
\frac{dp}{dt}=\frac{dp}{ds}\frac{ds}{dt}=-\frac{q}{||q||^{3}},\frac{dq}{ds}=%
\frac{dq}{ds}\frac{ds}{dt}=p.
\end{equation*}

 Since $\frac{(||p||^{2}+\varepsilon
h^{2})^{2}}{4}||q||^{2}=1$,
\begin{equation*}
E=\frac{1}{2}||p||^{2}-\frac{1}{||q||}=\frac{1}{2||q||}(||p||^{2}||q||-2)=%
\frac{1}{2||q||}(2-\varepsilon h^{2}||q||-2)=-\frac{1}{2}\varepsilon h^{2}.
\end{equation*}

So $E\,<0\ \ $\ in the spherical case and $E>0\ $in the hyperbolic case. 

The Euclidean case $E=0$ can be obtained by a limiting argument in which $%
\varepsilon \ $is regarded as a continuous parameter which tends to zero.\ \
To explain in more detail, let $\bar{x}(t)=x(t)-x_0(t)e_0$  where $x(t)$ is a solution of (\ref{geod1}).  If 

$w(t)=\lim_{\ep\rightarrow 0}\frac{1}{h^2\varepsilon }(\bar{x}(t))$, then $w(t)$ is a solution  of 
\[\frac{d^2 w}{dt^2}=0,
\]
that is, $w(t)$ is a geodesic corresponding to the standard Euclidean metric. It then follows from ({\ref{steregraphic}) that $\lim_{\ep\rightarrow 0}x_0=h$ and
 \[w=\lim_{\ep\rightarrow 0}\frac{1}{h^2\ep}\bar{x}=\lim_{\ep\rightarrow 0}\frac{2}{||p||^{2}+\varepsilon h^{2}%
}p=2\frac{p}{||p||^2}.\]  
Moreover,\ \ $\lim_{\varepsilon \rightarrow
0}dx_{0}=0\ $and \ $\lim_{\varepsilon \rightarrow 0}\frac{d\bar{x}}{\ep h^2}=dw=\frac{2}{||p||^{2}}%
dp-\frac{4p\cdot dp}{(||p||^{2})^{2}}p\ $ as can be seen from (\ref{diff1}\
).\ Therefore,\ $||dw||^{2}=\frac{4}{||p||^{4}}||dp||^{2}.$

The transformation $p\rightarrow w$ with $ w=\frac{2}{||p||^2}p\ $is the inversion about the circle $%
||p||^{2}=2\ $, and $||dw||^{2}=\frac{1}{||p||^{4}}||dp||^{2}$ is the
corresponding transformation of the Euclidean metric\ $||dp||^{2}$.\ The
Hamiltonian $H_0\ \ $associated with this metric is equal to
$\frac{1}{2}\frac{||p||^{4}}{4}||q||^{2}.$  

This Hamiltonian can be also obtained as the limit of $(\frac{h^2}{\ep})\frac{1}{2}\frac{(||p||^2+\ep h^2)^2}{4h^2}||q||^2$ when $\ep\rightarrow 0$.
On energy level $H=\frac{1}{%
2},\ \ ||p||^{2}||q||=2\ $ and therefore, $E=0.$

The integrals of motion for the problem of Kepler  are synonymous with the constancy of the matrix
$(x\wedge y)_{\varepsilon } $ along the flow of (\ref{wedge}). To be more specific, let 
 $x=x_{0}e_{0}+%
\bar{x}\ $ and  $y=y_{0}e_{0}+\bar{y}.$ Then,

$(x\wedge y)_{\varepsilon }=(x_{0}e_{0}+\bar{x}\wedge y_{0}e_{0}+\bar{y}%
)_{\varepsilon }=x_{0}(e_{0}\wedge \bar{y})_{\varepsilon }-y_{0(}e_{0}\wedge 
\bar{x})_{\varepsilon }+(\bar{x}\wedge \bar{y})_{\varepsilon }.$ \ Since ($%
x(t)\wedge y\ (t))_{\varepsilon }\ $is constant along the flows of (\ref{wedge})
both $x_{0}(e_{0}\wedge \bar{y})_{\varepsilon }-y_{0(}e_{0}\wedge \bar{x}%
)_{\varepsilon }$ \ and\ \ $(\bar{x}\wedge \bar{y})_{\varepsilon }\ $are \
also constant.\ But then the angular momentum
$L=q\wedge p\ $and the Runge-Lenz vector $F=Lp-\frac{q}{||q||}$ are given by 
\begin{equation}
L=(\bar{y}\wedge \bar{x})_{\varepsilon }\ and\ F=h(y_{0}(e_{0}\wedge \bar{x}%
)_{\varepsilon }-x_{0(}e_{0}\wedge \bar{y})_{\varepsilon })e_{0}.
\end{equation}

The first equality is evident from equations (\ \ref{steregraphic}) and (\ref%
{Hamlift}) and the fact that~ $(\bar{y}\wedge \bar{x})_{\varepsilon
}=\varepsilon (\bar{y}\wedge \bar{x}).$ Second equality follows by the
calculation below:

$(y_{0}(e_{0}\wedge \bar{x})_{\varepsilon }-x_{0(}e_{0}\wedge \bar{y}%
)_{\varepsilon })e_{0}=-y_{0}\bar{x}+x_{0}\bar{y}=$

$-(\frac{1}{h}q\cdot p)\frac{2\varepsilon h^{2}}{||p||^{2}+\varepsilon h^{2}}%
p+(\frac{h(||p||^{2}-\varepsilon h^{2})}{||p||^{2}+\varepsilon h^{2}})(\frac{%
||p||^{2}+\varepsilon h^{2}}{2h^{2}}q-\frac{q\cdot p}{h^{2}}p)=\frac{1}{h}%
(-(q\cdot p)p+\frac{||p||^{2}-\varepsilon h^{2}}{2}q).$

Since $\frac{(||p||^{2}+\varepsilon h^{2})}{2}||q||=1,\ \varepsilon h^{2}=%
\frac{2}{||q||}-||p||^{2},$ and therefore, 

$\ h(y_{0}(e_{0}\wedge \bar{x})_{\varepsilon
}-x_{0(}e_{0}\wedge \bar{y})_{\varepsilon })e_{0}=
(-(q\cdot p)p+||p||^{2}q-\frac{q}{||q||}=(q\wedge p)p-\frac{q}{||q||}=F.$

\subsubsection{Conic sections and the geodesics}
  The geodesics  of the spaces of constant curvature are
transformed into the conic sections of the problem of Kepler, a fact well known in conformal geometry. For the convenience of the reader not familiar with these facts and also for the completeness of the presentation we include the basic details. 

In the
spherical case, the great circle $x=a\cos \omega t+b\sin \omega t\ $with $%
||a||=||b||=h,\ a\cdot b=0$\  can be rotated around $e_{0}\ $\ so that $%
a $\ and $b\ $are in the \ subspace spanned by $e_{0},e_{1}, e_{2}.$
 Moreover, such a rotation $R$\ can be chosen so that $Ra=he_{1}.\ $
 
 Let $\alpha \ $%
denote the angle that the great circle makes with the plane $x_{0}=0.$Then,%
\begin{equation*}
x_{0}=h\sin \alpha \sin (ht),\ x_{1}=h\cos (ht),\ x_{2}=-h\cos \alpha \sin
(ht),\ x_{i}=0,\ i=3,\ldots ,n+1,
\end{equation*}

because $||x||^{2}||y||^{2}=h^{2}.$\ \ Furthermore, \ $\frac{h}{h-x_{0}}%
(x-he_{0})+he_{0}=(0,p)$\ \ implies that\ 
\begin{equation*}
p_{1}=\frac{h}{1-\sin \alpha \sin (ht)}\cos (ht),\ p_{2}=\frac{-h}{1-\sin
\alpha \sin (ht)}\cos \alpha \sin (ht),\ p_{i}=0,i=3,\ldots ,n,
\end{equation*}

and $ y(t)=\frac{1}{h^{2}}\frac{dx}{dt}$ implies that
\begin{equation*}
y=(\sin \alpha \cos (ht),-\sin (ht),-\cos
\alpha \cos (ht),\ 0,\ldots ,\ 0).
\end{equation*}

Then,%
\begin{equation*}
||p||^{2}+h^{2}=\frac{2h^{2}}{1-\sin \alpha \sin (ht)}\text{ and } \frac{p\cdot q}{h}=\sin \alpha
\cos (ht),
\end{equation*}
(implied by  equations\ (\ref{Hamlift})). Hence,%
\begin{equation*}\begin{array}{ll}&
y_{1}=\frac{q_{1}}{1-\sin \alpha \sin (ht)}-\frac{\sin \alpha \cos (ht)}{%
1-\sin \alpha \sin (ht)}\cos (ht),\\& y_{2}=\frac{q_{2}}{1-\sin \alpha \sin
(ht)}+\frac{\sin \alpha \cos (ht)}{1-\sin \alpha \sin (ht)}\cos \alpha \sin
(ht).
\end{array}
\end{equation*}

It follows that\ $-\sin (ht)\ (1-\sin \alpha \sin (ht))=q_{1}-\sin \alpha
\cos (ht)\cos (ht)$ \ and\ $-\cos \alpha \cos (ht)\ (1-\sin \alpha \sin
(ht))=q_{2}+\sin \alpha \cos (ht)\cos \alpha \sin (ht)$

 and therefore,%
\begin{equation*}
q_{1}=-\sin (ht)+\sin \alpha ,\ q_{2}=-\cos \alpha \cos (ht).
\end{equation*}

Hence, the great circle is transformed into the ellipse
\begin{equation*}
(q_{1}-\sin \alpha )^{2}+\frac{1}{\cos ^{2}\alpha }q_{2}^{2}=1
\end{equation*}

This ellipse  degenerates into a line \ through the origin when $%
\alpha =\frac{\pi }{2},$ or when the great circle passes through $he_{0}.$

A similar argument shows that the hyperboloid $x(t)=a\sinh (ht)+b\cosh (ht)\ 
$\ is transformed into  the hyperbola $-(q_{1}-\sin \alpha )^{2}+\frac{1}{\cos
^{2}\alpha }q_{2}^{2}=1.$

In the Euclidean case, the line $w=a+bt$ is transformed into\ the curve $\frac{2}{||p(t)||^2}p(t)$ via the mapping $w=\frac{2}{||p||^2}p$.  Hence,\ \ 

$b=\frac{dw}{dt}=\frac{2}{||p||^{2}}\frac{dp}{dt}-\frac{4}{||p||^{4}}%
(p, \frac{dp}{dt})p$.  After the substitutions, $(p,\frac{dp}{dt})=(w,\frac{dw}{dt}\frac{||p||^2}{||w||^2}$ and $\frac{dp}{dt}=\frac{||p||^4}{4}q$, $||p||^2=\frac{4}{||w||^2}$, the above equation becomes
\[b=2\frac{q}{||w||^2}+2\frac{1}{||w||^2}(w,\frac{dw}{dt})w.\]
Hence,
\[q=\frac{1}{2}(b(||a||^2-||b||^2t^2)-2a((a,b)+||b||^2t)).\]
This equation is a parabola in the $a,b$\ plane.

\section{Concluding remarks}

The above exposition \ could be viewed as a first step in unifying \ various
fragmented \ results in the theory of integrable systems. The fact that\
much of this theory is related to Lie groups and the associated Lie algebras
has been recognized in one form or another \ for some time now (\ \cite{Bol}%
\ \cite{per}\ \cite{rey}\ \cite{sem}\ ). \ However,\ in contrast to \ the
cited publications, the present study \ uses control theory and its Maximum
Principle as a point of departure for geometric problems with non-holonomic
constraints which greatly facilitates passage to the appropriate
Hamiltonians and which, at the same time, clarifies the role of the
Hamiltonians for the original problems.

Additionaly, the  ubiquitous presence of the affine problem on any symmetric space
paves a way for new classes of integrable systems, for it seems very likely
that the affine \ Hamiltonian is integrable on any coadjoint orbit ( \cite%
{Bol}).\ Further clarifications\ of this situation would be welcome
additions to the theory of integrable systems. \ Along more specific lines,
 the study of Fedorov and Jovanovic ((\cite%
{Fed}))  \ strongly\ suggests that the problem of Newmann on Steifel
manifolds \ can be seen also as the resticition of the affine Hamiltonian to
the coadjoint orbit of the semidirect product through an arbitrary symmetric
matrix.\ It would be instructive to investigate this situation in some
detail.

It might be also worthwhile to mention that the solutions of the affine
problem on the unitary group would find direct applications in the emerging
field of \ quantum control (\cite{Br}). \ This topic, however, because of
its own intricacies is deferred to a separate study.

\end{document}